\documentclass[preprint,11pt]{elsarticle}
\usepackage{amssymb,amsmath,amsthm,amsfonts,amscd}
\usepackage{graphicx,color}
\usepackage{cancel} 
\usepackage [latin1]{inputenc}
\journal{xxxxxx}
\newtheorem{theorem}{Theorem}[section]

\newtheorem{proposition}{Proposition}[section]

\newtheorem{remark}{Remark}[section]

\makeatletter
\def\tagform@#1{\maketag@@@{(\ignorespaces#1\unskip\@@italiccorr)}}
\makeatother

\newcommand{\cqdf}{\hfill\rule{5pt}{5pt}}
\newcommand{\jjnt}{\int\!\!\!\!\int}
\usepackage{fancyhdr}

\begin{document}

\begin{frontmatter}

\title{ Hierarchical Control  for the  Oldroyd Equation in Memoriam  to Professor Luiz Adauto Medeiros}

\author[Jesus]{Isa\'{\i}as Pereira  de Jesus \corref{cor1}}
\ead{isaias@ufpi.edu.br}
\author[Clark]{Marcondes Rodrigues Clark}
\author[Marinho]{Alexandro Marinho Oliveira }
\author[Trajano]{Aldo Trajano Lourêdo }
\address[Jesus]{Universidade Federal do Piau\'{\i}, DM, Campus Universitário Ministro Petrônio Portella, bairro Ininga, SG-4, Av. Universitária, s/n, 64049-550, Teresina, PI, Brazil, orcid: https://orcid.org/0000-0002-1917-1171}\cortext[cor1]{Corresponding author}
\address[Clark]{Universidade Federal do Piau\'{\i}, DM,  Campus Universitário  Ministro Petrônio Portella, bairro Ininga, SG-4, Avenida Universitária, s/n, 64049-550, Teresina, PI, Brazil,  orcid :
 https://orcid.org/0009-0009-4492-3628}
\address[Marinho]{Universidade Federal do Delta do Parna\'{\i}ba, DM, bairro Nossa Senhora de Fátima, Av. São Sebastião, 2819, 64202-020, Parna\'{\i}ba, PI, Brazil, orcid: https://orcid.org/0000-0003-3632-5133}
\address[Trajano]{Universidade Estadual da Para\'{\i}ba, DM, Rua Baraúnas, bairro Universitário, 351,  58429-500, Campina Grande, PB, Brazil, orcid https://orcid.org/0000-0002-3112-7456}


\begin{abstract}
This manuscript deals with a hierarchical control problem for Oldroyd equation under the Stackelberg-Nash strategy. The Oldroyd equation model is defined by non-regular coefficients, that is, they are
bounded measurable functions. We assume that we can act in the dynamic of the system by a hierarchy of controls, where one main control (the leader) and several additional secondary control (the followers) act in order to accomplish their given tasks: controllability for the leader and optimization for followers. We obtain the existence and uniqueness of Nash equilibrium and its characterization,  the approximate controllability  with respect
to the leader control, and the optimality system for leader control.
\end{abstract}

\begin{keyword}
Hierarchical Controllability; Stackelberg-Nash strategies; Oldroyd Fluid; Optimality System.
\MSC[2020] 93B05; 93C05; 35Q30; 76A05; 35Q93.\\

\end{keyword}

\end{frontmatter}
\section*{Introduction}
\label{intro}
\subsection*{Bibliographical comments}
With origin in game theory, and mainly motivated by economics, there exists several equilibrium concepts for multi-objective controllability of PDE. Each
of them determines a strategy. For example, we mention the noncooperative optimization strategy proposed by Nash (1951), the Pareto cooperative
strategy (1896),  and  the Von Stackelberg hierarchical-cooperative strategy (1934). The process in the problems above is a combination of strategies and is called {\it Stackelberg-Nash strategy.} The concept of {\it hierarchical control} was introduced by Lions (1994), where a simplified structure involving a single leader and a single follower control was considered to solve a problem of controllability for a hyperbolic equation.

For other equations, the hierarchical controllability has been considerably investigated. In the context of approximate controllability, we can cite some works involving Stackelberg-Nash strategy.  In fact, in the paper by  D\'iaz \& Lions (2005), the approximate distributed controllability of a parabolic system has been established following a Stackelberg-Nash strategy, and in L\'imaco et al. (2009), this same strategy was developed to obtain the approximate controllability for the  linear heat equation with moving boundaries.  In the context of linear fluid models, some investigations into approximate controllability using Stackelberg-Nash strategies began with the results of Guillén-González et al. (2013) for the Stokes system. These were later expanded upon by Araruna et al. (2014) for linearized micropolar fluids. Subsequently, Jesus \& Menezes (2015) extended the results of Araruna et al. (2014) to include moving domains. In both cases, the objective of the leader control is an approximate controllability result.

In recent years, much attention has been given to the investigation of  new
classes of problems in differential equations of hydrodynamics. Control problems
for the Navier-Stokes equations and other models of fluid mechanics are
examples of these . A considerable  number of papers and books (see, for
instance, Abergel \& Temam (1990), Betts  (2001), Blieger  (2007), Fursikov  (2000), Gunsburger (2003), Lions \& Zuazua (1998), Stavre  (2002),
and references therein) deal with the theoretical and numerical study of the above mentioned problems.  More precisely, in the context of fluid mechanics, significant controllability results are associated to the Burgers, Stokes, Euler, and Navier-Stokes equations. For instance, the local null controllability of the Burgers equation with distributed controls  was investigated in Fern\'andez-Cara \& Guerrero (2006), and recent work by Araruna et al. (2024) explores null control for Burgers equations within hierarchical controllability using the Stackelberg-Nash strategy.  For Stokes equations, the approximate and null controllability with distributed controls have been established in Fabre (1996) and  Imanuvilov (2001), respectively. Additionally, non-null controllability of Stokes equations with memory was analyzed by Fern\'andez-Cara et al. (2020). Global controllability results for Euler equations were proven by Coron (1996) and Glass (2000). Regarding Navier-Stokes equations, while only local exact controllability results are available for equations with initial and Dirichlet boundary conditions, as documented in Fursikov \& Imanuvilov (1999), Imanuvilov (2001), Fern\'andez-Cara et al. (2004, 2006), Coron et al. (2020) established a global exact controllability result for equations with Navier-slip (friction) boundary conditions.

The main novelties of this paper lie in the formulation and resolution of Nash and Stackelberg-Nash control problems within the framework of partial differential equations (PDEs) governing a non-conventional fluid model, specifically  the equations governing Oldroyd fluids . The results obtained here can generate several interesting problems generalizing or improving the results to other similar models, for instance,  this can be viewed as a first step in the path to understand similar questions in the context of the Navier-Stokes system and their variants.

\subsection*{The Oldroyd model}\label{sec11}
This model corresponds to an incompressible fluid  which is
described by the following system of partial differential equations:
\begin{equation}\label{eqI.1}
\begin{array}{c} \displaystyle \frac{\partial u}{\partial t} +
u.\nabla u +  \nabla p=div(\tau) + F(x,t),\,\,x\in \Omega,\,\,\,t>0,\vspace{0.3cm}\\
div(u)=0 ,\,\,x\in \Omega,\,\,\,t>0,
\end{array}
\end{equation}
with appropriate initial and boundary conditions. Here, $\Omega$ is an open
bounded connected set of $\mathbb{R}^n$ with smooth boundary $\Gamma$,
$\tau=(\tau_{ik})$ denotes the stress tensor with $tr\, \tau=0$, $u$
represents the velocity vector, $p$ is the pressure of the fluid and $F$ is an
external force. The stress tensor $\tau$ plays a special role because
the introduction of $\tau$ in $(\ref{eqI.1})$ has the purpose of letting us
consider reactions arising in the fluid during its motion. By establishing
(Hooke's Law) the connection between $\tau$, the tensor of deformation
velocities $D=(D_{ik})=\displaystyle\frac{1}{2}\left(u_{ix_k} + u_{kx_i}\right)$ and
their derivatives, we thus establish the type of fluid. Such relation between
$\tau$ and $D$ is called a \textit{defining or rheological equation or an
	equation of state} (see Serrin (1959)). From Newton's law, we set:
\begin{equation}\label{eqI.3}
\begin{array}{c}
\displaystyle \tau=2\nu D,
\end{array}
\end{equation}
where $\nu$ is the kinematic coefficient of viscosity. In this case, the fluid  is
called a  \textit{Newtonian Fluid}. Substituting $(\ref{eqI.3})$ into $(\ref{eqI.1})$ we
obtain the equations of motion of Newtonian fluids, which is called
\textit{Navier-Stokes equations}.

Over the last century and half, the model of a Newtonian fluid has been the basic
model of a viscous incompressible fluid. It describes flows of moderate velocities
of the majority of viscous incompressible fluids encountered in practice. However,
even earlier in the mid-nineteenth century it was known that there exists viscous
incompressible fluid not subject to the Newtonian equation $(\ref{eqI.3})$. That is,
it has a complex microstructure such as biological fluids, suspensions and liquid
crystals, which are used in the current industrial process and shows (nonlinear)
viscoelastic behavior that cannot be described by the classical linear viscous
Newtonian models. The first models of such fluids, were proposed in the
nineteenth century by Maxwell (1859, 1868), Kelvin(1875), Voigt (1889, 1892). In the mid-twentieth century, Oldroyd extended such models (see Oldroyd (1950, 1953, 1959, 1964)).

The model for  Oldroyd fluid (see Astarita \& Marruci  (1976), Wilkinson (1960))  can predict
the stress relaxation as well as the retardation of deformation. Due to this, it has
become popular for describing polymer suspension. To model the
behavior of a dilute polymer solution in a Newtonian solvent, the extra stress
tensor is often split into two components: a viscoelastic one and a purely viscous
one. So the Oldroyd fluids of order one as it is known in the Russian literature
(see Oskolkov (1989), Oskolkov \& Akhmatov (1974), Oskolkov \&  Kotsiolis (1986), Oskolkov et al. (1987)) are described by
defining relation:
\begin{equation}\label{eqI.4}
\begin{array}{c}
\displaystyle \left( 1 + \lambda \frac{\partial }{\partial t}
\right)\tau=2\nu \left(1 + k\nu^{-1}\frac{\partial }{\partial
	t}\right)D,
\end{array}
\end{equation}
where $\lambda,\,\nu,\,k$ are positive constants with  $\nu - k\lambda^{-1}>0$.
Here, $\nu$ denotes the kinematic viscosity, $\lambda$ is the relaxation time,
and $k$ represents the retardation time.

We observe that \eqref{eqI.4} can be rewritten in the form of an integral equation as follows:
\begin{equation}\label{eqI.5}
\begin{array}{c}
\displaystyle \tau(x,t)=2k\lambda^{-1}D(x,t) +
2\lambda^{-1}(\nu-k\lambda^{-1})
\int_0^te^{-\frac{(t-\sigma)}{\lambda}}D(x,\sigma)d\sigma,
\end{array}
\end{equation}
where $\tau(x,0) = 0.$

Thus, the equation of motion of the Oldroyd fluid of  first order can
be described most naturally by the system of integro-differential
equations:
\begin{equation}\label{eqI.6}
\begin{array}{l} \displaystyle \displaystyle \frac{\partial
u}{\partial t} + \left(u.\nabla\right) u - \mu \Delta u - \int_0^t
g(t-\sigma)\Delta u(x,\sigma)d\sigma + \nabla p= F(x,t) ,\,\,x\in \Omega,\,t>0, \vspace{0.3cm}
\end{array}
\end{equation}
and the incompressible condition:
\begin{equation}\label{eqI.7}
div(u)=0,\,\,x\in \Omega,\,t>0,
\end{equation}
with initial and boundary conditions:
\begin{equation}\label{eqI.8}
u(x,0)=u_0,\,\,x\in
\Omega,\,\,\,\,\mbox{and}\,\,\,\,u(x,t)=0\,\,x\in \Gamma,\,\,t\geq
0.
\end{equation}
Here  $\mu=k\lambda^{-1}$
and kernel $g(t)=\gamma e^{-\delta t}$, where
$\gamma=\lambda^{-1}(\nu - k\lambda^{-1}),$ with
$\delta=\lambda^{-1}$.  For details of the physical background and
its mathematical modeling, er refer to Oldroyd (1953), Astarita \& Marruci  (1976),
 Wilkinson (1960), and Oskolkov (1989).
\begin{remark}\label{isa557} As the theory of viscoelastic fluids describes flows with moderate velocities, the equation  $\eqref{eqI.6}$ admits a reasonable simplification, i.e., the convective term $\left(u.\nabla\right) u$ is neglected, as usual in mechanics.
\end{remark}
\subsection*{Notations}
As in Temam (1979)  let us denote $H^m(\Omega)$ as the standard
Hilbert-Sobolev space and by $\|\,.\,\|_m$ the norm defined on it. When $m=0$,
we call $H^0(\Omega)$ as the space of square-integrable functions
$L^2(\Omega)$ with the usual norm $|\,.\,|$ and inner product $(.,.)$. Further, let
$H^1_0(\Omega)$ be the completion of $C_0^{\infty}(\Omega)$ concerning
$H^1(\Omega)$-norm.

Let us consider the spaces:
\begin{equation*}
\begin{array}{l}
\mathcal{V}:=\left\{\phi\in
(C_0^{\infty}(\Omega))^n\,:\,div(\phi)=0\,\,\,\mbox{in}\,\,\,\Omega\right\},\\
H:=\,\mbox{the closure
	of}\,\,\mathcal{V}\,\,\mbox{in}\,\,(L^2(\Omega))^n-\mbox{space},
\end{array}
\end{equation*}
and
\begin{equation*}
\begin{array}{l}
V:=\,\mbox{the closure
	of}\,\,\mathcal{V}\,\,\mbox{in}\,\,(H^1_0(\Omega))^n-\mbox{space}.
\end{array}
\end{equation*}
The spaces of vector functions are indicated by boldface, for instance,
$\textbf{H}_0^1(\Omega)=(H_0^1(\Omega))^n$,
$\textbf{L}^2(\Omega)=(L^2(\Omega))^n$. The inner products on
$\textbf{H}_0^1(\Omega)$ and $\textbf{L}^2(\Omega)$ are defined by:
\begin{equation*}
\begin{array}{l}
((\phi,w)):=\sum_{i=1}^n(\nabla\phi_i,\nabla
w_i)\,\,\,\,\mbox{and}\,\,\,\,(\phi,w):=\sum_{i=1}^n(\phi_i, w_i),
\end{array}
\end{equation*}
respectively. Similarly , we define the norms:
\begin{equation*}
\begin{array}{l}
\|\phi\|:=\left(\sum_{i=1}^n|\nabla
\phi_i|^2\right)^{\frac{1}{2}}\,\,\,\,\mbox{and}\,\,\,\,|\phi|:=\left(\sum_{i=1}^n|\phi_i|^2\right)^{\frac12}.
\end{array}
\end{equation*}

Let us notice that under some smoothness assumptions on the boundary $\Gamma$, it
is possible to characterize $H$ and $V:$
\begin{equation*}
\begin{array}{l}
H=\{u\in\textbf{L}^2(\Omega); ~div (u)=0,\,\,\mbox{and}\,\,
u.\eta|_{\Gamma}=0\}\,\,\,\,\mbox{and}\,\,\,\,
V=\{u\in\textbf{H}_0^1(\Omega); ~div(u)=0\},
\end{array}
\end{equation*}
where $\eta = \eta(x)$ is the outward unit normal vector at $x \in \Gamma.$
By Poincar\'e's  inequality, it can be shown that the norm
$\textbf{H}_0^1(\Omega)$ is equivalent to
$\textbf{H}^1(\Omega)=(H^1(\Omega))^n-$norm. By $V'$ we denote the dual of
$V$.

\subsection*{Main Result}
Let $T>0$ be  a real number. We consider the cylindrical domain $Q := \Omega \times
(0,T)$ of $\mathbb{R}^{n+1}$ with lateral boundary $\sum := \Gamma \times
(0,T)$. We denote by $\mathcal{O},\mathcal{O}_1,\mathcal{O}_2,\ldots,\mathcal{O}_N$ non-empty
 disjoint  open subsets of $\Omega$.   By $\chi_{\mathcal{O}},$ $\chi_{\mathcal{O}_1},$ $\chi_{\mathcal{O}_2},$ \ldots, $\chi_{\mathcal{O}_N}$  we represent the characteristic
functions of $\mathcal{O}$, $\mathcal{O}_{1}$,  $\mathcal{O}_{2}$, \ldots, $\mathcal{O}_{N}$, respectively.

In this paper, we investigate the approximate controllability of the following system:
\begin{equation}\label{eqI.9}
\left|\begin{array}{l} \displaystyle \frac{\partial u}{\partial t} -
\mu\Delta u -
\displaystyle\int_0^tg(t-\sigma)\Delta u(\sigma)d\sigma+\nabla p=v\chi_{\mathcal{O}} +
\sum_{i=1}^Nw_i\chi_{\mathcal{O}_i}~~~~\mbox{in}~~~Q,\\
div(u)=0~~~~~\mbox{in}~~~~~Q,\\
u=0~~~~~\mbox{on}~~~~~\Sigma, \\
u(x,0)=u_0(x)~~~~\mbox{in}~~~~ \Omega,
\end{array}\right.
\end{equation}
where $u(x,t)=\left(u_1(x,t),\ldots,u_n(x,t)\right)$ is the  velocity  vector (or state of
the system) of moderate fluid evaluated at the point
$(x,t),\,x=(x_1,...,x_n)\in \mathbb{R}^n$, $p=p(x,t)$ is the pressure of the fluid
evaluated at the point $(x,t)$, $\mu$ represents a constant, and $u_0(x)$ is the initial
velocity.

The system $(\ref{eqI.9})$ can be interpreted as a variant of the classical Oldroyd equations $(\ref{eqI.6})-(\ref{eqI.8})$,  in which the nonlinearity $\left(u.\nabla\right) u$ has been omitted; see Remark $\ref{isa557}$.

In $(\ref{eqI.9}),$ the subset $\mathcal{O} \subset\Omega$ is the {\it main control domain}  (which is supposed to be as
small as desired), $\mathcal{O}_{1}$,  $\mathcal{O}_{2}$, \ldots, $\mathcal{O}_{N}$ are the {\it secondary control domains,}
the function $v$ is called  {\it leader control}, and $w_i$,  $(i=1,2, \ldots, N)$,  are the {\it followers controls.}

\begin{remark} By linearity of the system $(\ref{eqI.9}),$  without loss of generality, we may assume that $u_0 = 0.$
\end{remark}

As the solution $u$ of  $(\ref{eqI.9})$ depends on $v, w_1, \ldots, w_N$ then we denote it by
$\displaystyle u=u(x, t, v, \mathbf{w})$, where $\mathbf{w}=(w_1,\ldots, w_N)$, or
sometimes by $\displaystyle u=u(x, t, v, w_1,\ldots, w_N)$.


To localize the action of the controls $w_i$, $(i=1,2,\ldots, N)$, we introduce the functions
$\rho_i(x)$, defined in $\Omega$ with real values, satisfying:
\begin{equation}\label{eq4.2}
\begin{array}{c}
\rho_i\in L^{\infty}(\Omega),\,\,\rho_i\geq 0, \,\,\,\,
\rho_i=1\,\,\,\mbox{in} \,\,\,G_i\subset\Omega,
\end{array}
\end{equation}
where $G_i$ is a region where $w_i$ works.

We assume the leader objective to be of controllability type. On the other hand, the main objective of the followers is to hold $u(x,t,v,\textbf{w}),$ solution of the state equation  $(\ref{eqI.9})$
 at the time $T, $ near to a desired state $u^T(x),$ without a big cost for the controls $w_1,\ldots, w_N, $ with cost functionals defined  by:
\begin{equation}\label{eq4.3}
\begin{array}{c}
\displaystyle
J_i(v,\textbf{w}):=\frac12|w_i|^2_{\left(L^2(\mathcal{O}_i\times
	(0,T))\right)^n} +
\frac{\alpha_i}2\left|\rho_iu(.,T,v,\textbf{w})-\rho_iu^T(.)\right|^2_H,
\end{array}
\end{equation}
where $\alpha_i$ is a positive constant, $\textbf{w}=(w_1,....,w_N)$,
$v\in\left(L^2(\mathcal{O}\times (0,T))\right)^n$, and $w_i$ varies in
$\left(L^2(\mathcal{O}_i\times (0,T))\right)^n$,  with $i=1,2,\ldots,N.$

The Stackelberg-Nash strategy is described as follows: for each choice of the
leader $v$, we search for a Nash equilibrium for  the cost functionals $J_1,\ldots ,J_N,$
that is, we look for controls $w_1, . . . ,w_N$, depending
on $v,$ satisfying:
\begin{equation}\label{eq4.4}
\displaystyle J_i(v,w_1,\ldots,w_N)\leq
J_i(v,w_1,\ldots,\hat{w}_i,\ldots,w_N) \, \, \, \, \,\,\mbox{for
	all}\, \,\hat{w}_i\in \left(L^2(\mathcal{O}_i\times (0,T))\right)^n.
\end{equation}

The controls $w_1,\ldots, w_N$, solutions of the system of $N$ inequalities
$(\ref{eq4.4})$, are called \textit{Nash equilibrium} for the costs $J_1,\cdots,
J_N$ and they depend on $v$ (cf. Aubin (1984)).

For each $i = 1,2, \ldots, N, $ assuming  that $ \rho_i\in L^{\infty}(\Omega),$  and $\alpha_i$ is  small enough then there exists a unique Nash equilibrium $w_1,\ldots,w_N$, depending
on $v$, given by the inequalities $(\ref{eq4.4}).$   Later, we will make explicit this assumption with more details.

The main problems to be answered in this paper can be read as follows:

\textbf{$\bullet$ Problem 1} The existence of solutions $w_1,\ldots, w_N$ for the inequalities $(\ref{eq4.4}),$ that is, the existence of the Nash equilibrium for the functionals $J_1,\cdots,
J_N;$

\textbf{$\bullet$ Problem 2} Assuming that the existence of the Nash equilibrium $w_1(v), . . . ,w_n(v)$ was proved, then when $v$ varies in $\left(L^2(\mathcal{O}\times (0,T))\right)^n,$ to prove that the
solutions $u(x,t,v,\textbf{w}(v))$ of the state equation $(\ref{eqI.9})$, evaluated at $t = T$, that is, $u(x,T,v,\textbf{w}(v))$, generate a dense subset of $H$. This permits to approximate $u^T$.


Let us consider the following set:
\begin{equation*}
R(T):=\left\{\begin{array}{c} u(x,T,v,\textbf{w}(v)),\,v\in\left(L^2(\mathcal{O}\times
(0,T))\right)^n\,\,\mbox{where}\,\,
u(x,t,v,\textbf{w}(v))\,\,\mbox{is a}\\
\mbox{strong solution
	of}\,\,(\ref{eqI.9})\,\mbox{with}\,\,f=v\chi_{\mathcal{O}} +
\sum_{i=1}^Nw_i(v)\chi_{\mathcal{O}_i}\,\,\mbox{and}\,\,u_0=0.
\end{array}\right\}
\end{equation*}

Thus,  the main result of this paper is the following:
 \begin{theorem}\label{AC}
Let $ T > 0 $. Let us assume that for every $v\in \left(L^2(\mathcal{O}\times (0,T))\right)^n,$  there exists a unique Nash equilibrium $w_1,\ldots,w_N$, depending
on $v$, given by the inequalities $(\ref{eq4.4}),$ that is, for each $i = 1,2, \ldots, N, $ assume that $ \rho_i\in L^{\infty}(\Omega),$  and $\alpha_i$ is  small enough.  Then $R(T)$ is dense in $H$.

\end{theorem}

\subsection*{Organization of the paper}
This paper is organized as follows.  Initially, we present the optimality system for the followers controls. After this, we investigate the approximate controllability proving
the density Theorem \ref{AC}. The following section  is devoted to establish the existence and uniqueness of the Nash equilibrium by using the Lax-Milgram's Lemma.
Posteriorly,  we deal with the optimality system for the leader control. After this, we add some comments and point out open problems related to this article. Finally,  in the last section, we present the Appendix, where we prove the existence and uniqueness of solutions.

\section*{Optimal system for the followers controls}\label{sec2}
The main objective here is to express the followers controls $w_i$, $(i=1,2,\ldots, N)$, as weak
solutions of an appropriate system. For this, we suppose that there exists a Nash equilibrium $\textbf{w}=(w_1,\ldots,w_N)$ for the cost functionals  $J_i, (i=1,2,....,N),$ defined in $(\ref{eq4.3}).$
This means that $\textbf{w}=(w_1,\ldots,w_N)$ satisfies the Euler-Lagrange equation given by:
\begin{equation}\label{eq4.6}
J'_i(v,\textbf{w}(v)).\hat{w_i}=0\,\,\,\, \, \, \, \,\,\mbox{for
	all}\, \,\hat{w}_i\in \left(L^2(\mathcal{O}_i\times (0,T))\right)^n,
\end{equation}
where
$$J'_i(v,\textbf{w}(v)).\hat{w_i}=\frac{d}{d\lambda}
J_i\left(v,w_1,...,w_i+\lambda\hat{w}_i,...,w_N\right)\bigg|_{\lambda=0}.$$

The equation $(\ref{eq4.6})$ will allow us to reach to  the objective of this
section. In fact, from $(\ref{eq4.3}),$ we have that:

$\displaystyle J_i(v, w_1,..,w_i+\lambda\hat{w_i},\ldots,w_N)= \frac12\left(w_i
+\lambda\hat{w_i},w_i+\lambda\hat{w_i}\right)_{\left(L^2(\mathcal{O}_i\times(0,T))\right)^n}$
\begin{eqnarray*}
	&+&\displaystyle\frac{\alpha_i}{2}\bigg(\rho_iu(.,T,v,w_1,\ldots,w_i+\lambda\hat{w_i},\ldots,w_N) - \rho_iu^T(.),\rho_iu(.,T,v,w_1,\ldots,w_i +\lambda\hat{w_i},\ldots,w_N)\\
	&-& \rho_iu^T(.)\bigg)_H.
\end{eqnarray*}
As the state equation $(\ref{eqI.9})$ is linear, for any choice of the controls $ v, w_i$, its unique solution at the time $T$ can
be written as $  u(T) = L_0 v+ \sum_{i=1}^N L_i w_i,$ where $L_i$ are linear and continuous operators.  Later,
 we will make explicit this assumption with more details. Thus, the equality above can be rewritten as follows:\\
$\displaystyle
J_i(v,w_1,..,w_i+\lambda\hat{w_i},...,w_N)=\displaystyle\frac12\left(w_i+\lambda\hat{w_i},w_i+\lambda\hat{w_i}\right)_{\left(L^2(\mathcal{O}_i\times(0,T))\right)^n}$
\begin{equation}\label{eq4.7}
\begin{array}{rl}
+& \displaystyle\frac{\alpha_i}{2}\bigg(\rho_iL(w_i + \lambda\hat{w_i}) -\rho_iu^T(.), \rho_iL(w_i + \lambda\hat{w_i})
-\rho_iu^T(.)\bigg)_H,
\end{array}
\end{equation}
where $L(w_i + \lambda\hat{w_i}): = L_0v+L_1w_1+ ...+L_i( w_i+\lambda \widehat{w}_i)+ ...+L_Nw_N.$

Now, differentiating  the equation $(\ref{eq4.7})$ with respect to  $\lambda$  and
evaluate at
$\lambda=0,$ we obtain: \\
$\displaystyle \frac{d}{d\lambda}J_i(v,
w_1,\ldots,w_i+\lambda\hat{w_i},\ldots,w_N)\bigg|_{\lambda=0} \ = \
\left(w_i,\hat{w_i}\right)_{\left(L^2(\mathcal{O}_i\times(0,T))\right)^n}$
\begin{eqnarray*}
	&+& \displaystyle \alpha_i\bigg(\rho_i^2\left[u(.,T,v,w_1,\ldots,w_i,\ldots,w_N) -u^T(.)\right],L_i\hat{w}_i\bigg)_H,
\end{eqnarray*}
for all $\hat{w}_i \in L^2(\mathcal{O}_i\times (0,T))$, where
$L_i\hat{w_i}=\hat{u}_i(.,T,\hat{w_i})$ and $\hat{u}_i(x,t,\hat{w}_i)$ is the unique strong solution
of system:
\begin{equation}\label{eq4.8}
\left|\begin{array}{l} \displaystyle
\frac{\partial\hat{u}_i}{\partial t} - \mu\Delta \hat{u}_i -\int_0^tg(t-\sigma)\Delta \hat{u}_i(\sigma)d\sigma +
\frac{1}{N}\nabla p=\hat{w_i}\chi_{\mathcal{O}_i}~~\mbox{in}~~Q,\\
div(\hat{u}_i)=0~~~~~\mbox{in}~~~~~Q,\\
\hat{u}_i=0~~~~~\mbox{on}~~~~~\Sigma, \\
\hat{u}_i(x,0)=0~~~~\mbox{in}~~~~ \Omega,
\end{array}\right.
\end{equation}
associate to the  $\hat{w_i}$ with $N>0$ given.

Then, admitting  $(\ref{eq4.4})$ for the cost functionals $J_i, (i=1,2,....,N),$  we obtain the following
characterization:
\begin{equation}\label{eq4.9}
\begin{array}{c}
\displaystyle
\left(w_i,\hat{w_i}\right)_{\left(L^2(\mathcal{O}_i\times
	(0,T))\right)^n} + \displaystyle
\alpha_i\bigg(\rho_i^2\left[u(.,T,v,w_1,\ldots,w_i,\ldots,w_N) -
u^T(.)\right], \hat{u}_i(., T)\bigg)_H=0,
\end{array}
\end{equation}
for all $\hat{w}_i\in \left(L^2(\mathcal{O}_i\times (0,T))\right)^n.$

To find an optimality system for the followers, we represent by $\hat{p}=\hat{p}(x,t)\in L^2(0,T;L^2(\Omega))$ and
$\psi_i(x,t)=\left(\psi_{i1}(x,t),\ldots,\psi_{in}(x,t)\right)$, $(x,t)\in Q$ the weak solution of system:
\begin{equation}\label{eq4.10}
\left|\begin{array}{l} \displaystyle -\frac{\partial
	\psi_i}{\partial t} - \mu\Delta \psi_i-\int_t^Tg(\eta-t)\Delta\psi_i(\eta)d\eta + \nabla\hat{p}=0~~\mbox{in}~Q,\\
div(\psi_i)=0~~~~~\mbox{in}~~~~~Q, \\
\psi_i=0~~~~~\mbox{on}~~~~~\Sigma, \\
\psi_i(x,T)=\rho_i^2\left[u(x,T,v,\textbf{w}(v)) - u^T(x)\right]~~~~\mbox{in}~~~~
\Omega,
\end{array}\right.
\end{equation}
where the condition $\psi_i(x,T)$ is motivated by $(\ref{eq4.9})$.

Observe that if we make the change of variables $\tau=T-t$ in $(\ref{eq4.10})$ and set
$\varphi_i(x,\tau)=\psi_i(x,t),$ we transform $(\ref{eq4.10})$ into an equivalent
system in the  unknown $\varphi_i$ but with $\varphi_i(0)=\rho_i^2\left[u(.,T,v,\textbf{w}(v)) - u^T(.)\right]$,
where $\rho_i^2\left[u(.,T,v,\textbf{w}(v)) - u^T(.)\right] \in H.$ Furthermore, since the right hand side from $(\ref{eq4.10})_{1}$ is zero, then
$(\ref{eq4.10})$ admits a unique weak solution $\psi_i$ such that:
$$-\frac{\partial \psi_i}{\partial t} - \mu\Delta \psi_i -
\int_t^Tg(\eta-t)\Delta \psi_i(\eta)\,d\eta+ \nabla\widehat{p}=0,$$
in the sense of $L^2(0,T;V').$

If $\hat{u}_i(x,t,\hat{w}_i)$ is the strong solution of $(\ref{eq4.8})$,  we have that:
$$\hat{u}_i(.,.,\hat{w}_i)\in C^0([0,T];V)\cap L^2(0,T;V\cap
\textbf{H}^2(\Omega)).$$ Therefore, it makes sense to take the duality between
$$ \mathcal{L}(\psi_i)=-\frac{\partial \psi_i}{\partial t} -
\mu\Delta \psi_i - \int_t^Tg(\eta -t)\Delta \psi_i(\eta)\,d\eta +
\nabla\widehat{p}$$
and $\hat{u}_i:=\hat{u}_i(x,t,\hat{w}_i)$ in $L^2(0,T;V')\times L^2(0,T;V).$

Formally, we multiply both sides of $(\ref{eq4.10})$ by the strong solution
 $\hat{u}_i$ of $(\ref{eq4.8})$, and integrate in $Q$. Let us recall that we had assumed $u_0=0$ in $(\ref{eqI.9})$. Then:
\begin{eqnarray*}
	& & \displaystyle -\left(\psi_i(.,T),\hat{u}_i(T)\right)_H + \int_Q\psi_i \left\{\frac{\partial \hat{u}_i}{\partial t} - \mu\Delta \hat{u}_i\right\}dxdt +\displaystyle \int_{0}^T\left<\nabla\widehat{p}(t),\hat{u}_i(t)\right>_{(H^{-1}(\Omega))^n\times (H_0^1(\Omega))^n}dt \\
	& &-\int_0^T\left<\int_t^Tg(\eta-t)\Delta\psi_i(\eta)d\eta,\hat{u}_i(t)\right>_{V'\times V}dtdt=0.
\end{eqnarray*}
From $(\ref{eq4.8}),$ it follows that:
$$
\frac{\partial \hat{u}_i}{\partial t} -\mu\Delta \hat{u}_i =
\hat{w}_i\chi_{\mathcal{O}_i}-\frac{1}{N}\nabla p +
\int_0^tg(t-\sigma)\Delta \hat{u}_i(\sigma)d\sigma \;\; \,\,\mbox{a.e.
	in}\;\; \,\,Q.
$$
Hence,
\begin{equation}\label{eq4.11}
\begin{array}{ll}
-&\displaystyle\bigg(\rho_i^2\left[\hat{u}_i(.,T,v,\textbf{w}(v)) -u^T(.)\right],\hat{u}_i(.,T)\bigg)_H + \int_Q\psi_i\hat{w}_i\chi_{\mathcal{O}_i}dxdt \\
+&\displaystyle\int_{0}^T\left<\nabla\widehat{p}(t),\hat{u}_i(t)\right>_{(H^{-1}(\Omega))^n\times (H_0^1(\Omega))^n}dt \\
+& \displaystyle\int_0^T\left(\int_0^tg(t-\sigma)\Delta\hat{u}_i(\sigma)\,d\sigma,\psi_i(t)\right)dt \\
-& \displaystyle\int_0^T\left<\int_t^Tg(\eta-t)\Delta\psi_i(\eta)d\eta,\hat{u}_i(t)\right>_{V'\times V}dt\\
-& \displaystyle\frac1N\int_0^T\left(\nabla p(t),\psi_i(t)\right)dt=0.
\end{array}
\end{equation}

\begin{remark}\label{obs2.3}
	Let us notice that:
	\begin{equation*}
	\begin{array}{l}
	\displaystyle \left<\nabla p(t),u(t)\right>_{(H^{-1}(\Omega))^n\times (H_0^1(\Omega))^n} = -\left(p(t),div u(t)\right)_{L^2(\Omega)}=0,
	\end{array}
	\end{equation*}
	since $u(t)\in V$ and thus $div (u(t)) = 0.$
\end{remark}
\begin{remark}\label{obs2.4}
	By Fubini's formula,  we have that:
	\begin{equation*}
	\displaystyle\int_0^T\left(\int_0^tg(t-\sigma)\Delta
	u(\sigma)\,d\sigma,\psi(t)\right)dt-\int_0^T\left<\int_t^Tg(\eta-t)\Delta
	\psi(\eta)d\eta,u(t)\right>_{V'\times V}dt=0.
	\end{equation*}
\end{remark}

From Remarks $\ref{obs2.3}$ and  $\ref{obs2.4}$, we have that the equation
$(\ref{eq4.11})$ becomes:
\begin{equation*}
\begin{array}{l}
\displaystyle \bigg (\rho_i^2\left[u(.,T,v,\textbf{w}(v)) -
u^T(.)\right], \hat{u}_i(., T)\bigg)_H = \int_Q\psi_i
\hat{w}_i\chi_{\mathcal{O}_i} dxdt,\;  \mbox{ for all } \hat{w}_i\in\left(L^2(\mathcal{O}_i\times (0,T))\right)^n.
\end{array}
\end{equation*}
Hence, from $(\ref{eq4.9})$, we deduce that:
\begin{equation*}
\begin{array}{l}
\left(w_i + \alpha_i\psi_i,\hat{w_i}\right)_{\left(L^2(\mathcal{O}_i\times
(0,T))\right)^n} =0,\;\;\  \; \mbox{ for all } \hat{w}_i\in\left(L^2(\mathcal{O}_i\times (0,T))\right)^n,
\end{array}
\end{equation*}
i.e.,
\begin{equation}\label{eq4.13}
w_i=-\alpha_i\psi_i\,\, a.e.\,\,\mbox{in}\,\,\,\,\mathcal{O}_i\times(0,T).
\end{equation}
\begin{remark}
	The relation $(\ref{eq4.13})$ is very  important. It corresponds to finding the control
	functions as a weak solution of  system $(\ref{eq4.10}).$ The
	main reason to express the followers controls in such way is to find an optimality system. Using this fact, one can derive numerical approximation algorithms. 
\end{remark}
Thus, the best $w_i$ satisfying $(\ref{eq4.4})$ are given by $w_i=-\alpha_i\psi_i,$
where $\psi_i$ is the unique solution of the following optimality system:
\begin{equation}\label{eq4.14}
\left|\begin{array}{l} \displaystyle \frac{\partial u}{\partial t} -
\mu\Delta u - \int_0^tg(t-\sigma)\Delta u(\sigma)d\sigma +\nabla p + \sum_{i=1}^N\alpha_i\psi_i\chi_{\mathcal{O}_i}=
v\chi_{\mathcal{O}}~\mbox{in}~Q,\vspace{0.3cm}\\
\displaystyle -\frac{\partial\psi_i}{\partial t} - \mu\Delta \psi_i  - \int_t^Tg(\eta-t)\Delta \psi_i(\eta)d\eta + \nabla\widehat{p}=0~~~~\mbox{in}~~~Q,\vspace{0.3cm}\\
div(u)=0,\,\,div(\psi_i)=0~~~~~\mbox{in}~~~~~Q, \vspace{0.3cm}\\
u=0,\,\,\,\psi_i=0~~~~~\mbox{on}~~~~~\Sigma, \vspace{0.3cm}\\
u(x,0)=0,\,\,\,\,\, \psi_i(x,T)=\rho_i^2\left[u(x,T,v,\textbf{w}(v))
- u^T(x)\right]~~~~\mbox{in}~~~~
\Omega.
\end{array}\right.
\end{equation}

\section*{On approximate controllability}\label{sec3}
Our goal in this section is to show the approximate controllability for the state
equation $(\ref{eqI.9})$ assuming that  $v\in \left(L^2(\mathcal{O}\times
(0,T))\right)^n$, i.e.,  we will prove that the solution $u(x,t,v,\textbf{w}(v))$ of
the state problem $(\ref{eqI.9}),$ evaluated at $t = T,$ generate a  dense subset of  $H$. Observing the
cost functionals  $J_i, (i=1,2,....,N),$ defined in $(\ref{eq4.3}),$ we will have the
approximate controllability with $u^T\in H$ as in  D\'iaz \& Lions (2005).

Now, we present the proof of Theorem \ref{AC}.

\noindent{\bf Proof of Theorem \ref{AC} }  By linearity of optimality system $(\ref{eq4.14}),$  without loss of generality, let us assume that $u^T\equiv0$ (it suffices to use a
translation argument). We will prove that if $f\in H$ satisfies
$\left(u(.,T,v,\textbf{w}(v)),f\right)_H=0$ for all $v\in\left(L^2(\mathcal{O}\times
(0,T))\right)^n$, then $f$ is the  null vector of $H$, i.e., the orthogonal
complement of $R(T)$ in $H$ is null.

In fact, multiplying $(\ref{eq4.14})_1$ by $\varphi,$ $(\ref{eq4.14})_2$ by
$\xi_i,$ respectively, and  integrating in $Q$, we obtain:
\begin{equation}\label{eq4.15}
\begin{array}{ll}
& \displaystyle \left(u(.,T,v,\textbf{w}(v)),\varphi(T)\right)_H + \int_Q u\left\{-\frac{\partial \varphi}{\partial t} - \mu\Delta\varphi \right\}dxdt
+\displaystyle \int_0^T(\nabla p(t),\varphi(t))dt\vspace{0.2cm}\\
- &\displaystyle\int_0^T\left(\int_0^tg(t-\sigma)\Delta u(\sigma)\,d\sigma,\varphi(t)\right)dt
+\displaystyle\int_Q\varphi\sum_{i=1}^N\alpha_i\psi_i\chi_{\mathcal{O}_i}dxdt \vspace{0.2cm}\\
= &\displaystyle\int_Q\varphi v\chi_{\mathcal{O}}dxdt,
\end{array}
\end{equation}
and
\begin{equation}\label{eq4.16}
\begin{array}{ll}
& \displaystyle  \int_Q \psi_i\left\{\frac{\partial \xi_i}{\partial t}- \mu\Delta \xi_i\right\}dxdt + \left(\psi_i(0),\xi_i(0)\right)_H
+\int_0^T\left<\nabla\hat{p}(t),\xi_i\right>dt \vspace{0.2cm}\\
-&\displaystyle\int_0^T\left<\int_t^Tg(\eta-t)\Delta\psi_i(\eta)d\eta,\xi_i(t)\right>dt -\left(\psi_i(T),\xi_i(T)\right)_H =0,
\end{array}
\end{equation}
where $i=1,2,\ldots,N.$ Fixing the conditions:
\begin{equation}\label{eq4.17}
\left|\begin{array}{l} \displaystyle -\frac{\partial
	\varphi}{\partial t} - \mu\Delta \varphi- \int_t^Tg(\eta-t)\Delta
\varphi(\eta)d\eta + \nabla\hat{p}
=0~~~~\mbox{in}~~~Q,\vspace{0.3cm}\\
\displaystyle \frac{\partial
	\xi_i}{\partial t} - \mu\Delta \xi_i - \int_0^tg(t-\sigma)\Delta
\xi_i(\sigma)d\sigma +  \nabla {p}=-\alpha_i\varphi\chi_{\mathcal{O}_i}~~\mbox{in}~~~Q,\\
div(\varphi)=0,\,\,div(\xi_i) =0~~~~~\mbox{in}~~~~~Q, \\
\varphi=0,\,\,\,\xi_i=0~~~~~\mbox{on}~~~~~\Sigma, \\
\xi_i(x,0)=0,\,\,\,\,\, \varphi(x,T)=f(x) +
 \sum_{i=1}^N\rho_i^2\left[\xi_i(x,T)\right]~~~~\mbox{in}~~~~
\Omega,
\end{array}\right.
\end{equation}
we conclude that:

\begin{itemize}
	\item[(I)] From $(\ref{eq4.17})_1$, the equation $(\ref{eq4.15})$ can be rewritten as follows:
	\begin{equation}\label{eq4.18}
	\begin{array}{ll}
	& \displaystyle \left(u(.,T,v,\textbf{w}(v)),\varphi(T)\right)_H - \int_0^T\left(\int_0^tg(t-\sigma)\Delta u(\sigma)\,d\sigma,\varphi(t)\right)dt \vspace{0.2cm} \\
	+ & \displaystyle \int_0^T\left(\nabla p(t),\varphi(t)\right)dt - \int_0^T\left < u(t), \nabla\hat{p}(t) \right>dt  + \int_Q\sum_{i=1}^N\varphi\alpha_i\psi_i\chi_{\mathcal{O}_i}dxdt\vspace{0.2cm}\\
	+ & \displaystyle\int_0^T\left<\int_t^Tg(\eta-t)\Delta\varphi(\eta)d\eta,u(t)\right>dt =\int_Q\varphi v\chi_{\mathcal{O}}dxdt.
	\end{array}
	\end{equation}
	\item[(II)] From $(\ref{eq4.17})_2$, the expression in  $(\ref{eq4.16})$ becomes:
	\begin{equation}\label{eq4.19}
	\begin{array}{ll}
	- &\displaystyle\left(\psi_i(T),\xi_i(T)\right)_H - \int_Q \psi_i\alpha_i\varphi\chi_{\mathcal{O}_i}dxdt  + \int_0^T\left<\nabla\hat{p}(t),\xi_i\right>dt \vspace{0.2cm}\\
	- & \displaystyle\int_0^T\left<\int_t^Tg(\eta-t)\Delta\psi_i(\eta)d\eta,\xi_i(t)\right>dt - \int_0^T\left(\psi_i,\nabla p(t)\right)dt \vspace{0.2cm}\\
	+ & \displaystyle\int_0^T\left(\int_0^tg(t-\sigma)\Delta \xi_i(\sigma)d\sigma,\psi_i(t)\right)dt=0.
	\end{array}
	\end{equation}
\end{itemize}

\begin{remark} The system $(\ref{eq4.17})$ admits a unique solution $\{\varphi,\xi_i,\widehat{p},p\}$ in $C^0([0,T];H)\times
	C^0([0,T];V)\times L^2(0,T;L^2(\Omega))\times  L^2(0,T;H^1(\Omega))$;  see Appendix. Therefore, it makes sense the calculus in $(\ref{eq4.15})$ and $(\ref{eq4.16})$.
\end{remark}

From $(\ref{eq4.14})_3,$ $(\ref{eq4.17})_3,$ $(\ref{eq4.18})$, $(\ref{eq4.19}),$ and Remarks  $\ref{obs2.3}$ and $\ref{obs2.4},$  we have that:
\begin{itemize}
	\item[(I)]  Equation $(\ref{eq4.18})$ is reduced to the expression:
	\begin{equation}\label{eq4.20}
	\begin{array}{c}
	\displaystyle \left(u(.,T,v,\textbf{w}(v)),\varphi(T)\right)_H +
	\int_Q\sum_{i=1}^N\varphi\alpha_i\psi_i\chi_{\mathcal{O}_i}dxdt
	=\int_Q\varphi v\chi_{\mathcal{O}}dxdt.
	\end{array}
	\end{equation}
	\item[(II)] Expression in  $(\ref{eq4.19})$ becomes:
	\begin{equation}\label{eq4.21}
	\begin{array}{c}
	\displaystyle  - \left(\psi_i(T),\xi_i(T)\right)_H - \int_Q
	\psi_i\alpha_i\varphi\chi_{\mathcal{O}_i}dxdt=0.
	\end{array}
	\end{equation}
\end{itemize}

Now, summing $(\ref{eq4.21})$ from $1$ at $N$, we obtain that:
\begin{equation}\label{eq4.22}
\begin{array}{c}
\displaystyle  \sum_{i=1}^N\left(-\psi_i(T),\xi_i(T)\right)_H
=\int_Q \sum_{i=1}^N\psi_i\alpha_i\chi_{\mathcal{O}_i}\varphi dxdt.
\end{array}
\end{equation}

Hence, substituting $(\ref{eq4.22})$ in $(\ref{eq4.20})$, we get:
\begin{equation}\label{eq4.23}
\begin{array}{c}
\displaystyle \left(u(.,T,v,\textbf{w}(v)),\varphi(T)\right)_H +
\sum_{i=1}^N\left(-\psi_i(T),\xi_i(T)\right)_H=\int_Q\varphi
v\chi_{\mathcal{O}}dxdt.
\end{array}
\end{equation}

From $(\ref{eq4.14})_5$ and $(\ref{eq4.17})_5,$ we find that:
\begin{equation}\label{eq4.200}
\begin{array}{c}
\varphi(T)=f +
\displaystyle \sum_{i=1}^N\rho_i^2\xi_i(T),\,\,\,\,\,\psi_i(T)=\rho_i^2u(T).
\end{array}
\end{equation}

Thus, combining   $(\ref{eq4.23})$ and $(\ref{eq4.200})$, we deduce that:
$$
\begin{array}{c}
\displaystyle(u(T),f)_H +
\left(u(T),\sum_{i=1}^N\rho_i^2\xi_i(T)\right)_H -
\sum_{i=1}^N\left(\rho^2_i u(T),\xi_i(T)\right)_H=\int_Q\varphi
v\chi_{\mathcal{O}}dxdt,
\end{array}
$$

i.e,
\begin{equation}\label{eq.100}
\begin{array}{c}
\displaystyle(u(T),f)_H=\int_Q\varphi v\chi_{\mathcal{O}}dxdt,\,\,\,\mbox{for all}\,\,\, v\in\left(L^2(\mathcal{O}\times (0,T))\right)^n.
\end{array}
\end{equation}

From assumption  $(u(T),f)_H=0$,  then from $(\ref{eq.100})$, it follows that:
$$
\begin{array}{c}
\displaystyle\int_Q\varphi v\chi_{\mathcal{O}}dxdt=0,\,\,\,\mbox{for all}\,\,\, v\in\left(L^2(\mathcal{O}\times (0,T))\right)^n ,
\end{array}
$$
and therefore,
\begin{equation}\label{eq.1000}
\begin{array}{l}
\varphi=0\,\,\,\,\,\mbox{in}\,\,\,\,\mathcal{O}\times (0,T)\subset Q.
\end{array}
\end{equation}
It follows by the  unique continuation (cf. Doubova  \& Fern\'andez-Cara  (2012), Lemma 2.3, p. 576), that  $\varphi = 0 \; in\; Q$. Thus, coming back to $(\ref{eq4.17})$, we
obtain that $\xi_i = 0 \; in\; Q,$ $\,\,\forall\,i=1,2,\ldots,N. $ Therefore, as
$\varphi\in C^0([0,T];H)$ and $\xi_i\in C^0([0,T];V)$,  we have that $\xi_i(x,T)=0$
and $\varphi(x,T)=0,$ and from $(\ref{eq4.17})_5$ implies that $f\equiv0$ in
$\Omega$. This ends the proof of Theorem \ref{AC}. \cqdf

\section*{Existence and uniqueness of  Nash equilibrium}\label{sec4}
Our goal in this section is to prove the existence of Nash equilibrium  for the cost functionals  $J_i, (i=1,2,\ldots,N),$ defined in $(\ref{eq4.3}),$
corresponding to the state equation $(\ref{eqI.9})$. For this, let
$H_i=\left(L^2(\mathcal{O}_i\times (0,T))\right)^n$ and $
\mathcal{H}=\prod_{i=1}^N H_i$ be Hilbert spaces. Moreover, we consider for
each $i=1,...,N,$ the cost functionals $J_i,$ given by $(\ref{eq4.3}),$ and the operator $
L_i:\left(L^2(\mathcal{O}_i\times (0,T))\right)^n \to V $ such that
$L_i\hat{w}_i=\hat{u}_i(T)$, where $\hat{u}_i$ is a unique strong solution of the linear
system $(\ref{eq4.8}).$ The functionals $L_i$ are linear and are
well defined. Indeed, since $\hat{u}_i$ is the strong solution of the system $(\ref{eq4.8})$
with $\hat{w}_i$, it follows that $ \hat{u}_i$ belongs to $C^0([0,T];V)$. Hence, $\hat{u}_i(T)\in
V$. Moreover, $ \|\hat{u}_i\|_{C^0([0,T];V)}\leq C_T\|\hat{w}_i\|_{H_i},$ and thus, $
\|\hat{u}_i(T)\|_V\leq C_T\|\hat{w}_i\|_{H_i}, $ that is, $ \|L_i\hat{w}_i\|_V\leq
C_T\|\hat{w}_i\|_{H_i}, $ for all $i=1,\ldots,N$. Therefore $L_i\in \mathcal{L}(H_i,V)$.
Since $V\subset H$ with continuous and dense immersion, we have that
$|L_i\hat{w}_i|_H\leq C_0\|L_i\hat{w}_i\|_V\leq \tilde{C}_0\|\hat{w}_i\|_{H_i},$
where $\tilde{C}_0=C_0C_T.$

Let  $u$ be a strong solution of the problem $(\ref{eqI.9})$ with $u_0=0.$ From
$(\ref{eq4.8})$,  we have:
\begin{equation}\label{eq4.25}
\left|\begin{array}{l} \displaystyle \sum_{i=1}^N\frac{\partial u_i}{\partial t} - \sum_{i=1}^N\mu\Delta u_i
-\displaystyle\int_0^tg(t-\sigma)\sum_{i=1}^N \Delta u_i(\sigma)d\sigma
+\displaystyle\sum_{i=1}^N\frac{1}{N}\nabla p \vspace{0.2cm}\\ = \displaystyle\sum_{i=1}^Nw_i(v)\chi_{\mathcal{O}_i}~~~~\mbox{in}~~~Q, \vspace{0.2cm}\\
\displaystyle\sum_{i=1}^Ndiv(u_i)=0~~~~~\mbox{in}~~~~~Q,\vspace{0.2cm}\\
\displaystyle\sum_{i=1}^Nu_i=0~~~~~\mbox{on}~~~~~\Sigma, \vspace{0.2cm}\\
\displaystyle\sum_{i=1}^Nu_i(x,0)=0~~~~\mbox{in}~~~~ \Omega.
\end{array}\right.
\end{equation}
Since $v\in \left(L^2(\mathcal{O}\times (0,T))\right)^n, $ we compare $(\ref{eqI.9})$
with $(\ref{eq4.25})$ and conclude that $ u=z + \sum_{i=1}^Nu_i, $ where $u_i$
is the strong solution of $(\ref{eq4.8})$   with $w_i=\hat{w}_i$ and $z$ is fixed
depending only of $v$.

Since $u\in C^0([0,T];V)$, we have that $ u(x,T,v,\textbf{w}(v))=z^T +
\sum_{i=1}^N L_iw_i, $ for $z^T$ fixed.

This notation allows us  to rewrite the cost functionals  $J_i, (i=1,2,\ldots,N),$ defined in $(\ref{eq4.3}),$ in the
form:
\begin{equation}\label{eq4.26}
\begin{array}{c}
\displaystyle J_i(v,\textbf{w})=\frac12|w_i|^2_{H_i} +
\frac{\alpha_i}2\left|\rho_i\left(\sum_{j=1}^NL_jw_j
-\eta^T(.)\right)\right|^2_H,
\end{array}
\end{equation}
where $\eta^T=z^T-u^T.$

This means that $\textbf{w}= (w_1,\ldots,w_N) \in \mathcal{H}$ will  be a Nash
equilibrium for the convex functionals $J_i$ given by $(\ref{eq4.3})$, if  their
Gateaux derivatives are null in any direction $\hat{w_i}$, i.e., we must show
that $\textbf{w}\in \mathcal{H}$ satisfies:
\begin{equation*}
\begin{array}{l}
\displaystyle \left(w_i,\hat{w_i}\right)_{H_i} + \displaystyle
\alpha_i\left(\rho_i\left[\sum_{j=1}^NL_jw_j
-\eta^T(.)\right],\rho_iL_i\hat{w}_i\right)_H=0, \;\;\ \mbox {for all}\; \hat{w}_i\in H_i,
\end{array}
\end{equation*}
or
\begin{equation*}
\begin{array}{l}
\displaystyle \left(w_i,\hat{w_i}\right)_{H_i} + \displaystyle
\alpha_i\left(L_i^{*}\left[\rho_i^2\sum_{j=1}^NL_jw_j\right],\hat{w}_i\right)_H
- \alpha_i\left(L_i^{*}\left(\rho_i^2
\eta^T\right),\hat{w}_i\right)_H=0
\end{array}
\end{equation*}
for all $\hat{w}_i\in H_i,$ where $L^{*}_i\in \mathcal{L}(H,H_i)$ is the  adjoint
operator of the  $L_i\in \mathcal{L}(H_i,H)$. Hence,
\begin{equation}\label{eq4.27}
\begin{array}{c}
\displaystyle w_i +
L_i^{*}\left[\rho_i^2\sum_{j=1}^NL_jw_j\right]=\alpha_iL_i^{*}\left(\rho_i^2
\eta^T\right)\,\,\,\mbox{in}\,\,\,H_i,\,\,\,\mbox{for
	all}\;\;\;\;\ \,i=1,2,\ldots,N.
\end{array}
\end{equation}

Since $\alpha_iL_i^{*}\left(\rho_i^2 \eta^T\right)\in H_i$, we can  aim to find a vector
$$\left( \alpha_1L_1^{*}\left(\rho_1^2\eta^T\right),\ldots,\alpha_NL_N^{*}\left(\rho_N^2
\eta^T\right) \right)$$ in $\mathcal{H}$. After that, we define
the functional $A:\mathcal{H}\to \mathcal{H}$ such that $ \left(Aw,
\hat{w}\right)_{\mathcal{H}}=\left(f,\hat{w}\right)_{\mathcal{H}},\,\,\,\forall\,f\in
\mathcal{H},\,\,\forall\,\hat{w}\in \mathcal{H}$, where $A\in
\mathcal{L}(\mathcal{H},\mathcal{H})$ is given by $ (Aw)_i=w_i +
L_i^{*}\left[\rho_i^2\sum_{j=1}^NL_jw_j\right]. $ Thus, the problem can be
formulated as the following question:

\textbf{$\bullet$ Question :} Given $f\in \mathcal{H}$, is there a unique
$w\in\mathcal{H}$  such that $\left(Aw,\hat{w}\right)_{\mathcal{H}}=
\left(f,\hat{w}\right)_{\mathcal{H}},\,\,\forall\,\,\hat{w}\in \mathcal{H}$ ?

Therefore, we must prove that the linear equation
$\left(Aw,\hat{w}\right)_{\mathcal{H}}= \left(f,\hat{w}\right)_{\mathcal{H}}$ admits a
solution $w=(w_1,\cdots,w_N)$ in $\mathcal{H}$ for each $f=(f_1,....,f_N)$ in
$\mathcal{H}$. To guarantee  its solvability, we will apply  the Lax-Milgram's
Lemma, with certain restrictions on $\alpha_i$ and $\rho_i$. More precisely, the following result holds:
	\begin{proposition}\label{is600}
Let us assume  that
\begin{equation*}
\begin{array}{l}
\rho_i\in L^{\infty}(\Omega),\;\;\alpha_i = \alpha \;\; \mbox{for  i = 1,\ldots, N,}\;\;{\mbox and}\;\;\alpha\;\;{\mbox is}\;\;{\mbox small}\;\;{\mbox enough}
\end{array}
\end{equation*}
such that $\beta_0:=\tilde{C}_0^2\alpha\max_{i,j=1,\ldots,N}\|\rho_i-\rho_j\|_{L^{\infty}(\Omega)}\max_{i=1,\ldots,N}\|\rho_i\|_{L^{\infty}(\Omega)}<1.$
Then there exists a Nash equilibrium $\textbf{w}= (w_1,\ldots, w_N)$ for the functionals $J_i$ defined in $(\ref{eq4.3})$.
\end{proposition}
\proof
	Let us observe that:
	\begin{eqnarray*}
		\displaystyle\left(Aw,w\right)_{\mathcal{H}} &=& \sum_{i=1}^N|w_i|_{H_i}^2 + \displaystyle\sum_{i=1}^N\alpha_i\left(\sum_{j=1}^N\rho_jL_jw_j,\rho_iL_iw_i\right)_H \\
		& & + \sum_{i,j=1}^N\alpha_i\left((\rho_i-\rho_j)L_jw_j,\rho_iL_iw_i\right)_H.
	\end{eqnarray*}
	Therefore, according to the hypothesis $\alpha=\alpha_i,\,\,\mbox{for
	all}\;\; \,i=1,2,\ldots,N$, we
	obtain:
	\begin{equation}\label{eq4.28}
	\begin{array}{lll}
	\displaystyle\left(Aw,w\right)_{\mathcal{H}} &=&  \displaystyle \alpha\left|\sum_{i=1}^N\rho_iL_iw_i\right|_H^2 +\sum_{i=1}^N|w_i|_{H_i}^2 \\
	& & +\displaystyle \alpha\sum_{i,j=1}^N\left((\rho_i-\rho_j)L_jw_j,\rho_iL_iw_i\right)_H.
	\end{array}
	\end{equation}
	Next, let us that:
	\begin{equation*}
	\begin{array}{l}
	\displaystyle\alpha\left|\sum_{i,j=1}^N\left((\rho_i-\rho_j)L_jw_j,\rho_iL_iw_i\right)_H\right|_{\mathbb{R}}\leq
	\displaystyle\tilde{C}_0^2\alpha\max_{i,j=1,\ldots,N}\|\rho_i-\rho_j\|_{L^{\infty}(\Omega)}
	\max_{i=1,..,N}\|\rho_i\|_{L^{\infty}(\Omega)}|w|_{\mathcal{H}}^2.
	\end{array}
	\end{equation*}
	Hence, from $(\ref{eq4.28})$, we obtain:
	\begin{equation}\label{eq4.29}
	\begin{array}{l}
	\displaystyle\left(Aw,w\right)_{\mathcal{H}} \geq (1-\beta_0)|w|_{\mathcal{H}}^2,
	\end{array}
	\end{equation}
	where $\displaystyle
	\beta_0=\tilde{C}_0^2\alpha\max_{i,j=1,\ldots,N}\|\rho_i-\rho_j\|_{L^{\infty}(\Omega)}
	\max_{i=1,..,N}\|\rho_i\|_{L^{\infty}(\Omega)}.$
	Since $\alpha$ is small enough, we have $\gamma_0=(1-\beta_0)>0$. Therefore,
	from $(\ref{eq4.29})$, it follows that  $ \left(Aw,w\right)_{\mathcal{H}}\geq
	\gamma_0|w|^2_{\mathcal{H}}.$
	
	Finally, since $A\in \mathcal{L}(\mathcal{H},\mathcal{H})$ and
	$\left(Aw,w\right)_{\mathcal{H}}\geq \gamma_0|w|^2_{\mathcal{H}}$, we have from
	Lax-Milgram's Theorem that, for a given $f\in \mathcal{H}$, there exists
	a unique $w\in \mathcal{H}$ such that $Aw=f$. In particular, for $f=\left(
	\alpha_1L_1^{*}\left(\rho_1^2\eta^T\right),\ldots,\alpha_NL_N^{*}\left(\rho_N^2
	\eta^T\right) \right)\in \mathcal{H}$, there exists a unique solution $w\in
	\mathcal{H}$, that is, a Nash equilibrium for the cost functionals $J_i,$  satisfying $
	(Aw)_i=\alpha_iL_i^{*}\left(\rho_i^2\eta^T\right),\; \mbox{for
	all}\;\;\,i=1,2,\ldots,N$, i.e.,
  \begin{equation*}
	\begin{array}{l}
	\displaystyle w_i +
	L_i^{*}\left[\rho_i^2\sum_{j=1}^NL_jw_j\right]=\alpha_iL_i^{*}\left(\rho_i^2
	\eta^T\right)\,\,\,\mbox{in}\,\,\,H_i, \;\; \,\mbox{for
	all}\;\; \,i=1,2,\ldots,N.
	\end{array}
	\end{equation*} \cqdf

\section*{Optimal system for the leader control} \label{sec5}
In the previous sections we have seen  that no matter what  strategy  the leader
assumes, the followers make their choices $w_1,....,w_N$  satisfying  the Nash
equilibrium. Moreover, this choices can be made through of an optimality
system. The goal of this section is to obtain an optimality system for the leader control.
For this, we consider the functional: $$ J(v)=\frac12\int_{\mathcal{O}\times (0,T)} v^2\, dx\,dt,
$$
and the minimization problem
\begin{equation}\label{eq4.30}
\,\,\left\{
\begin{array}{l}
\inf J(v)\\
\mbox{subject to}\,\,\,u(T,v)\in u^T + \varepsilon B,
\end{array}\right.
\end{equation}
where $\varepsilon>0$ is a given real number, $B\equiv B(0,1)$ is
the unitary ball of $H,$ and $u(v)$ is the unique solution of the
optimality system defined in $(\ref{eq4.14}).$

We introduce two convex proper functions as follows:
\begin{itemize}
	\item{} The first one is defined in $\left(L^2(\mathcal{O}\times (0,T))\right)^n$
	by:
	\begin{equation}\label{eq4.31}
	\displaystyle F_1(v):=\frac12\int_{\mathcal{O}\times (0,T)} v^2\,
	dx\,dt.
	\end{equation}
	\item{} The second one is  defined in $H$ by:
	\begin{equation}\label{eq4.32}
	F_2(f):=\left\{
	\begin{array}{l} 0\,\,\,\,if\,\,\,f\in u^T + \varepsilon
	B,\\
	+ \infty \,\,\,\,out.
	\end{array}\right.
	\end{equation}
\end{itemize}
With these notations, problem $(\ref{eq4.30})$ is equivalent to
\begin{equation}\label{eq4.33}
\displaystyle\inf_{v \in \left(L^2(\mathcal{O}\times (0,T))\right)^n}\left\{F_1(v) + F_2(L(v))\right\},
\end{equation}
where $L(v)=u(.,T,v)$, with  $L\in \mathcal{L}\left(\left[L^2(\mathcal{O}\times
(0,T))\right]^n, H\right)$.

Observe that there exists  $v\in \left(L^2(\mathcal{O}\times (0,T))\right)^n$
such that $F_1$ and $F_2$ are finite,  $F_1$ is continuous in $v,$ and  $F_2$ is continuous in $L(v)$. By the Duality Theorem of Fenchel and Rockafellar (1969)  (see also Brezis (2010), Ekeland \&  Temam (1974)), we have:
\begin{equation}\label{eq4.34}
\displaystyle \inf_{v \in \left(L^2(\mathcal{O}\times (0,T))\right)^n}\left(F_1(v)+
F_2(Lv)\right)=-\inf_{f \in H}\left(F_1^*(L^*f) + F_2^*(-f)\right),
\end{equation}
where $L^*$ denotes the adjoint of $L$ and $F_i^*$ is the conjugate
function of $F_i\;\; (i = 1,2).$ This means that the
primal problem (\ref{eq4.30}) is equivalent to its dual.

Now, our next goal is to get the best control function of the leader. In fact, using
$(\ref{eq.100})$ then  for $f\in H,$  it follows that:
\begin{equation*}
\begin{array}{l}
\displaystyle(L(v),f)_H=\int_Q\varphi v\chi_{\mathcal{O}}dxdt,   \,\,\, \mbox{for all}\,\,\ v \in\left(L^2(\mathcal{O}\times (0,T))\right)^n.
\end{array}
\end{equation*}
Thus, we conclude that:
\begin{equation}\label{eq4.35}
L^*f=\varphi\chi_{\mathcal{O}},
\end{equation}
where $\varphi$ is a solution of $(\ref{eq4.17}).$

We see easily that:
\begin{equation}\label{eq4.36}
F_1^*(v)=F_1(v),
\end{equation}
and
\begin{equation}\label{eq4.37}
F_2^*(f)=\left(f,u^T\right)_H + \varepsilon \|f\|,
\end{equation}
where $\|.\|$ denotes the  norm of $H$.

 So, from (\ref{eq4.35}) -- (\ref{eq4.37}), the expression in (\ref{eq4.34}) becomes:
\begin{equation}\label{eq4.38}
\displaystyle \inf_{v \in \left(L^2(\mathcal{O}\times (0,T))\right)^n}\left(F_1(v)+
F_2(Lv)\right)=-\inf_{f \in H} F(f),
\end{equation}

where the functional $F: H \rightarrow \mathbb{R}$ is defined by:
\begin{equation}\label{isa5509}
\begin{array}{l}
\displaystyle F(f):=\frac12 \int_{\mathcal{O}\times (0,T)}(L^*f)^2\,dx\,dt +
\varepsilon\|f\|-\left(f,u^T\right)_H,
\end{array}
\end{equation}
with $L^*f$  given by (\ref{eq4.35}).

Let us notice that the problem (\ref{eq4.30})  has a unique solution (see Rockafellar (1967)). Consequently, the dual problem  has also a unique solution.

After a quick computation of the Gateaux derivative of the functional $(\ref{isa5509})$, then for $\hat{\varphi}\; \in
H,$  we obtain the following variational inequality (cf.  Ekeland \&  Temam (1974)):
\begin{equation}\label{eq4.40}
\begin{array}{l}
\displaystyle\int_{\mathcal{O}\times (0,T)}\varphi(\hat{\varphi} -
\varphi)\,dx\,dt\, + \varepsilon\|\hat{f}\| - \varepsilon\|f\|-
\left(\hat{f}-f,u^T\right)_H\geq 0,  \quad\forall \hat{f} \in
H.
\end{array}
\end{equation}
Now, for each $f\in H$, we consider the unique solution of $(\ref{eq4.17})$ $\varphi$ and introduce $(u,\psi_i)$, $i=1,2,\ldots,N$, as the unique
solutions of the system:
\begin{equation}\label{eq4.41}
\left|\begin{array}{l} \displaystyle \frac{\partial u}{\partial t} -
\mu\Delta u - \int_0^tg(t-\sigma)\Delta u(\sigma)d\sigma +\nabla p
+ \sum_{i=1}^N\alpha_i\psi_i\chi_{\mathcal{O}_i}=
\varphi\chi_{\mathcal{O}}~\mbox{in}~Q,\vspace{0.3cm}\\
\displaystyle -\frac{\partial \psi_i}{\partial t} - \mu\Delta \psi_i
- \int_t^Tg(\eta-t)\Delta
\psi_i(\eta)d\eta \,+\nabla\widehat{p}=0~~~~\mbox{in}~~~Q,\vspace{0.3cm}\\
div(u)=0,\,\,div(\psi_i)=0~~~~~\mbox{in}~~~~~Q, \vspace{0.3cm}\\
u=0,\,\,\,\psi_i=0~~~~~\mbox{on}~~~~~\Sigma, \vspace{0.3cm}\\
u(x,0)=0,\,\,\,\,\,
\psi_i(x,T)=\rho_i^2\left[u(x,T,v)\right]~~~~\mbox{in}~~~~
\Omega.
\end{array}\right.
\end{equation}
Multiplying $(\ref{eq4.41})_1$ by $\hat{\varphi}-\varphi$, and $(\ref{eq4.41})_2$ by
$\hat{\xi}_i -\xi_i,$ respectively, and  integrating in $Q$, we obtain:
\begin{equation}\label{eq4.42}
\begin{array}{ll}
& \displaystyle \left(u(.,T),\hat{\varphi}(T)-\varphi(T)\right)_H + \int_0^T\left(\nabla p(t),\left(\hat{\varphi}-\varphi\right)\right)dt\,\\
+ & \displaystyle \int_Q u\left\{-\frac{\partial\left(\hat{\varphi}-\varphi\right)}{\partial t} -\mu\Delta\left(\hat{\varphi}-\varphi\right)\right\}dxdt \\
- & \displaystyle \int_0^T\left<u(t),\int_t^Tg(\eta-t)\Delta\left(\hat{\varphi}-\varphi\right)(\eta)\,d\eta\right>dt \, \\
+ & \displaystyle\int_Q\left(\hat{\varphi}-\varphi\right)\sum_{i=1}^N\alpha_i\psi_i\chi_{\mathcal{O}_i}dxdt \\
= & \displaystyle\int_Q\varphi\left(\hat{\varphi}-\varphi\right)\chi_{\mathcal{O}}dxdt,
\end{array}
\end{equation}
and
\begin{equation}\label{eq4.43}
\begin{array}{ll}
& \displaystyle\left(\psi_i(0),\left(\hat{\xi}_i(0)-\xi_i(0)\right)\right)_H +\int_0^T\left<\nabla\hat{p}(t),\left(\hat{\xi}_i-\xi_i\right)\right>dt \\
+ & \displaystyle  \int_Q \psi_i\left\{\frac{\partial\left(\hat{\xi}_i-\xi_i\right)}{\partial t} - \mu\Delta \left(\hat{\xi}_i-\xi_i\right)\right\}dxdt \\
- & \displaystyle\int_0^T\left(\psi_i(t),\int_0^tg(t-\sigma)\Delta \left(\hat{\xi}_i-\xi_i\right)(\sigma)d\sigma\right)dt \\
+ & \left(\psi_i(T),\hat{\xi}_i(T)-\xi_i(T)\right)_H =0,
\end{array}
\end{equation}
where $i=1,2,\ldots,N.$

Setting $\bar{\varphi}=\hat{\varphi}-\varphi$, $\bar{\xi}_i=\hat{\xi}_i-\xi$, and fixing the
conditions:
\begin{equation}\label{eq4.44}
\left|\begin{array}{l} \displaystyle -\frac{\partial
	\bar{\varphi}}{\partial t} - \mu\Delta \bar{\varphi} -
\int_t^Tg(\eta-t)\Delta \bar{\varphi}(\eta)d\eta
=0~~~~\mbox{in}~~~Q,\vspace{0.3cm}\\
\displaystyle \frac{\partial
	\bar{\xi}_i}{\partial t} - \mu\Delta \bar{\xi}_i - \int_0^tg(t-\sigma)\Delta
\bar{\xi}_i(\sigma)d\sigma=-\alpha_i\bar{\varphi}\chi_{\mathcal{O}_i}~~\mbox{in}~~~Q,\vspace{0.3cm}\\
div(\bar{\varphi})=0,\,\,div(\bar{\xi}_i)=0~~~~~\mbox{in}~~~~~Q, \vspace{0.3cm}\\
\bar{\varphi}=0,\,\,\,\bar{\xi}_i=0~~~~~\mbox{on}~~~~~\Sigma, \vspace{0.3cm}\\
\bar{\xi}_i(x,0)=0,\,\,\,\,\, \bar{\varphi}(x,T)=\hat{f}(x) - f(x) +
 \displaystyle \sum_{i=1}^N\rho_i^2\left[\bar{\xi}_i(x,T)\right]~~~~\mbox{in}~~~~
\Omega,
\end{array}\right.
\end{equation}
we deduce that:
\begin{equation}\label{eq4.45}
\begin{array}{ll}
& \displaystyle \left(u(.,T),\hat{\varphi}(T)-\varphi(T)\right)_H + \int_0^T\left(\nabla p(t),\hat{\varphi}(t)-\varphi(t)\right)dt \\
+ & \displaystyle\int_Q\sum_{i=1}^N\left(\hat{\varphi}(t)-\varphi(t)\right)\alpha_i\psi_i\chi_{\mathcal{O}_i}dxdt
=\int_Q\varphi \left(\hat{\varphi}-\varphi\right)\chi_{\mathcal{O}}dxdt
\end{array}
\end{equation}
and
\begin{equation}\label{eq4.46}
\begin{array}{ll}
- & \displaystyle\left(\psi_i(T),\hat{\xi}_i(T)-\xi_i(T)\right)_H - \int_Q\psi_i\alpha_i\left(\hat{\varphi}(t)-\varphi(t)\right)\chi_{\mathcal{O}_i}dxdt\\
+ & \displaystyle\int_0^T\left<\nabla\hat{p}(t),\hat{\xi}_i-\xi_i\right>dt=0,
\end{array}
\end{equation}
where $i=1,2,\ldots,N.$
From Remark $\ref{obs2.3}$,  we have that the expressions in  $(\ref{eq4.45})$ and $(\ref{eq4.46})$ becomes:
\begin{equation}\label{eq4.47}
\begin{array}{ll}
& \displaystyle \left(u(.,T),\hat{\varphi}(T)-\varphi(T)\right)_H + \int_Q\sum_{i=1}^N\left(\hat{\varphi}(t) - \varphi(t)\right) \alpha_i\psi_i\chi_{\mathcal{O}_i}dxdt \\
= & \displaystyle\int_Q\varphi \left(\hat{\varphi}-\varphi\right)\chi_{\mathcal{O}}dxdt
\end{array}
\end{equation}
and
\begin{equation}\label{eq4.48}
\begin{array}{l}
\displaystyle  - \left(\psi_i(T),\hat{\xi}_i(T)-\xi_i(T)\right)_H -
\int_Q
\psi_i\alpha_i\left(\hat{\varphi}(t)-\varphi(t)\right)\chi_{\mathcal{O}_i}dxdt=0,\,\,\,
\end{array}
\end{equation}
where $i=1,2,\ldots,N.$\\
Summing $(\ref{eq4.48})$ from $1$ to $N$ and substituting into $(\ref{eq4.47})$,
we get:
\begin{equation}\label{eq4.49}
\begin{array}{ll}
 \displaystyle \left(u(.,T),\hat{\varphi}(T)-\varphi(T)\right)_H + \sum_{i=1}^N\left(-\psi_i(T),\hat{\xi}_i(T)-\xi_i(T)\right)_H
= \displaystyle\int_Q\varphi\left(\hat{\varphi}-\varphi\right)\chi_{\mathcal{O}}dxdt.
\end{array}
\end{equation}
From $(\ref{eq4.41})_5$ and $(\ref{eq4.44})_5$, it follows that:
\begin{equation}\label{eq4.400}
\begin{array}{ll}
\hat{\varphi}(T)-\varphi(T)=\hat{f}-f +
\sum_{i=1}^N\rho_i^2\left(\hat{\xi}_i(T)-\xi_i(T)\right),\,\,\,\,\,\psi_i(T)=\rho_i^2u(.,T).
\end{array}
\end{equation}
Combining $(\ref{eq4.49})$ and $(\ref{eq4.400})$, we find that:
\begin{equation}\label{eq4.50}
\left(u(.,T),\hat{f}-f\right)_H=\int_{\mathcal{O}\times(0,T)}\varphi\left(\hat{\varphi
}-\varphi\right)\,dxdt.
\end{equation}
Substituting $(\ref{eq4.50})$ in $(\ref{eq4.40})$,  we get:
\begin{equation}\label{eq4.51}
\left(u(.,T)-u^T,\hat{f}-f\right)_H + \varepsilon\|\hat{f}\|
-\varepsilon\|f\|\geq 0,\,\,\,\,\mbox{for all}\,\,\hat{f}\in H.
\end{equation}
More precisely,  we summarize these results in the following theorem:
\begin{theorem}\label{teo5.1}
	The best control function $v$ of the leader, that is, the function that minimizes
	$$
	\frac12\int_{\mathcal{O}\times (0,T)}v^2\,dx\,dt
	$$
	subject to $u(.,T,v,\textbf{w}(v))\in u^T + \varepsilon B,$ is given by:
	$$
	v=\varphi\chi_{\mathcal{O}},
	$$
	where $\varphi$ is given from unique solution $\left\{u,\psi_i,\varphi,\xi_i\right\}$ of
	the optimality system:
	$$
	\left|\begin{array}{l} \displaystyle \frac{\partial u}{\partial t} -
	\mu\Delta u - \int_0^tg(t-\sigma)\Delta u(\sigma)d\sigma +\nabla p + \sum_{i=1}^N\alpha_i\psi_i\chi_{\mathcal{O}_i}=
	\varphi\chi_{\mathcal{O}}~\mbox{in}~Q,\vspace{0.3cm}\\
	\displaystyle -\frac{\partial
		\psi_i}{\partial t} - \mu\Delta \psi_i  - \int_t^Tg(\eta-t)\Delta \psi_i(\eta)d\eta +
	\nabla\widehat{p}=0~~~~\mbox{in}~~~Q,\vspace{0.3cm}\\
	\displaystyle -\frac{\partial
		\varphi}{\partial t} - \mu\Delta \varphi - \int_t^Tg(\eta-t)\Delta
	\varphi(\eta)d\eta +  \nabla\hat{p}
	=0~~~~\mbox{in}~~~Q,\vspace{0.3cm}\\
	\displaystyle \frac{\partial
		\xi_i}{\partial t} - \mu\Delta \xi_i - \int_0^tg(t-\sigma)\Delta
	\xi_i(\sigma)d\sigma +  \nabla{p}=-\alpha_i\varphi\chi_{\mathcal{O}_i}~~\mbox{in}~~~Q,\vspace{0.3cm}\\
	div(\varphi)=0,\,\,div(\xi_i)=0~~~~~\mbox{in}~~~~~Q, \vspace{0.3cm}\\
	div(u)=0,\,\,div(\psi_i)=0~~~~~\mbox{in}~~~~~Q, \vspace{0.3cm}\\
	\varphi=0,\,\,\,\xi_i=0,\,\,\,u=0,\,\,\,\psi_i=0~~~~~\mbox{on}~~~~~\Sigma, \vspace{0.3cm}\\
	\xi_i(x,0)=0,\,\,\,\,\, \varphi(x,T)=f(x) +
	\displaystyle \sum_{i=1}^N\rho_i^2\left[\xi_i(x,T)\right]~~~~\mbox{in}~~~~
	\Omega,\vspace{0.3cm}\\
	u(x,0)=0,\,\,\,\,\,
	\psi_i(x,T)=\rho_i^2\left[u(x,T,v)\right]~~~~\mbox{in}~~~~
	\Omega,
	\end{array}\right.
	$$
	\vspace{0.3cm} with $f$ being the unique solution of the  variational inequality:
	\begin{equation}\label{eq4.555}
	\left(u(.,T,f)-u^T,\hat{f}-f\right)_H + \varepsilon\|\hat{f}\|
	-\varepsilon\|f\|\geq 0,\,\,\,\,\mbox{for all}\,\,\hat{f}\in H.
	\end{equation}
\end{theorem}

\begin{remark}
	In $(\ref{eq4.555})$ we rewrite $u(.,T,f)$ to emphasize the fact that the solutions
	$\left\{u,\psi_i,\varphi,\xi_i\right\}$ of the optimality system also depends on $f.$
\end{remark}
\section*{Conclusions}\label{sec6}
 In this  this section we make some comments and briefly discuss some possible extensions of our results and also indicate open issues on the subject.

\subsection*{Extensions to various analogous scenarios.}
We have seen in this paper that it is possible to obtain approximate controllability result for the Oldroyd equation following a Stackelberg-Nash strategy. The point of adding secondary controls $w_1,..., w_N$ consists
  in the task of being a low cost control that ensures that the solution u of (\ref{eqI.9}) is not far from  ideal state  $u^T.$ Importantly, the concepts and techniques employed  in this paper can be applied to various
  analogous scenarios, such as hierarchical  controllability for linear and semilinear parabolic and hyperbolic equations, different types of non-cylindrical control domains, similar boundary control problems, among others. Indeed, concerning these themes, several prior studies are noteworthy. For example, in Límaco et al. (2009) the authors present the  hierarchical approximate controllability for the linear heat equation in a non-cylindrical domain. Jesus (2015) demonstrated hierarchical approximate controllability for linearized micropolar fluids in moving domains, employing a Stackelberg-Nash strategy. In the case of an Oldroyd fluid system, the exact controllability of Galerkin's approximations is proposed and analyzed in Marinho et al. (2014). In   Araruna et al. (2015), the authors  established a Stackelberg-Nash strategy with exact controllability for the leader control in  semilinear parabolic equations, utilizing Carleman inequalities. For semilinear hyperbolic equations, a Stackelberg-Nash strategy with exact controllability for the leader control
is proved in Araruna et al. (2018).

\subsection*{Extensions to Navier-Stokes equations.}
 The controllability  of Navier-Stokes equations have received significant attention in recent years, as evidenced by the extensive study documented in Coron et al. (2020) along with the latest advancements. However, to our knowledge, there exists a gap in the literature regarding null controllability for these equations. For further insights, see Araruna et al. (2015), p. 20, or more recent discussions in Araruna et al. (2024), p. 15. This unresolved issue presents an intriguing and pertinent avenue for future research.

\subsection*{Extension to  nonlinear Oldroyd system.}
    The theoretical analysis of nonlinear Oldroyd systems has been a focal point of research efforts in recent years, as evidenced by the works of Lions \& Masmoudi (2000), Fern\'andez-Cara et al. (2002),
    Galdi (2008), and Renardy (2009). Notably, Fern\'andez-Cara et al. (2020), p. 4, emphasize the significance of investigating Oldroyd fluids governed by nonlinear partial differential equations.
    Consequently, it is interesting to explore how the findings of this paper might be applicable in the context of the nonlinear Oldroyd fluids system, for instance, in the context of the equations  (\ref{eqI.6})-- (\ref{eqI.8}). Note that this system in question is  more difficult to solve than the system (\ref{eqI.9}). The main reason is the presence of the nonlinear term $\left(u.\nabla\right) u$ in (\ref{eqI.6}).  However, whether or not the results in this paper can be extended to this framework is at present an open question.

\subsection*{The case with a different definition for cost functionals (\ref{eq4.3}).}
	In  Guillén-González et al. (2013), the authors solve an approximate control problem for the Stokes equation using the Stackelberg-Nash strategy. In this article, the cost functionals for the followers controls  $v^{i},\;  (i = 1,2),$ are defined by:
	\begin{equation}
		\label{funcmedar}
		J_i\left(f,\,v^1,\,v^2\right)
		\,:=\,
		\dfrac{\alpha_i}{2}\jjnt_{\omega_{i,d}\times(0,T)}|u-u_{i,d}|^2\,dx\,dt
		+
		\dfrac{\mu_i}{2}\jjnt_{\omega_i\times(0,T)}|v^i|^2\,dx\,dt,
	\end{equation}
	where $\alpha_i,\,\mu_i$ are positive constants, and  $u_{i,d}$ are functions given in $L^2(\omega_{i,d}\times(0,T)).$
	
	Thus, a priori, the leader and the followers have different tasks, with the followers's task being to prevent the state function from deviating from a given function.

Note that the situation considered in this article is different: here, all the  controls have the same objective. However, we assume a situation where depending on the region, the control may need to have a different configuration, that is, if we have two regions, we will need two controls. Therefore, similar to the article by Guillén-González et al. (2013), we assume a hierarchy among the controls and employ with the Stackelberg optimization (cooperative) strategy, with a leader and followers. In the case of three or more regions, we will have one leader and $N$ followers ($N>1$), and they will act without collaboration among themselves, hence we use the Nash strategy.
	
Moreover, if we follow the same procedures as in this article but with the functionals defined as in \eqref{funcmedar}, we will have the followers controls  characterized by the following optimality system:
	\begin{equation}
		\label{otimedar}
		\left|
		\begin{array}{lll}
			\displaystyle \dfrac{\partial u}{\partial t}-\mu\Delta u - \int_0^tg(t - \sigma)\Delta u(\sigma)\,d\sigma + \nabla p
			\,=\,
			v\chi_{\mathcal{O}} - \sum_{i=1}^{2}\dfrac{1}{\mu_i}\,q^i\,\chi_{\omega_i}
			& \mbox{in} & Q,\\
            \displaystyle-\dfrac{\partial q^i}{\partial t}-\mu\Delta q^i - \int_t^Tg (\sigma - t)\Delta q^i(\sigma)\,d\sigma + \nabla\widehat{p}
			\,=\,
			\alpha_{i}(u-u_{i,d})\chi_{\omega_{i,d}}  & \mbox{in}& Q, \\
            \displaystyle  div(u)=0,\,\,div(q^i)=0    &\mbox{in}& Q, \\
			\displaystyle u\,=\,0,\quad
			q^{i} = 0 &  \mbox{on} & \Sigma, \\
			\displaystyle u(x\,,0)\,=\, u_{0},\quad q^{i}(x\,,T)=0 & \mbox{in} & \Omega.\\
		\end{array}
		\right.
	\end{equation}
Finally, we can solve the approximate control problem for the system \eqref{otimedar} following the same ideas as in this article.

\subsection*{Non null controllability for the system  (\ref{eqI.9}).}
It is expectable that the system (\ref{eqI.9}) (with a nonzero initial data for $u$) is not  null controllable, in view of other previous results on the controllability of parabolic systems with memory;  see  Guerrero \& Imanuvilov (2013)  for more details. Thus, it would also be quite interesting to obtain a result asserting that null controllability for the system \eqref{eqI.9} does not hold in general, which we plan to present
in a forthcoming paper.

On the other hand, exact controllability problems also can be considered in this context.  In general, exact controllability does not hold. Another interesting issue is what happens as the viscosity coefficient goes to zero. It is known that, for $\mu = 0$, exact controllability holds at least under some geometric control conditions. For interested readers on this subject, we cite for instance Boldrini et al. (2012).


\section*{Appendix: On existence and uniqueness of solutions}\label{sec7}
This appendix aims to prove the existence and uniqueness of solutions to the
coupled system  $(\ref{eq4.17})$. To simplify the notation, we will make the change of variable
$\phi(x,t)=\varphi(x,T-t)$ into $(\ref{eq4.17}),$ which yields the following equivalent system:
\begin{equation}\label{eqA.1}
\left|\begin{array}{l} \displaystyle \frac{\partial
	\phi}{\partial t} - \mu\Delta \phi \;\ -\int_0^tg(t-\eta)\Delta
\phi(\eta)d\eta + \nabla\hat{p}
=0~~~~\mbox{in}~~~Q,\vspace{0.3cm}\\
\displaystyle \frac{\partial
	\xi_i}{\partial t} - \mu\Delta \xi_i - \int_0^tg(t-\sigma)\Delta
\xi_i(\sigma)d\sigma +  \nabla {p}=-\alpha_i\varphi\chi_{\mathcal{O}_i}~~\mbox{in}~~~Q,\\
div(\phi)=0,\,\,div(\xi_i) =0~~~~~\mbox{in}~~~~~Q, \\
\phi=0,\,\,\,\xi_i=0~~~~~\mbox{on}~~~~~\Sigma, \\
\xi_i(x,0)=0,\,\,\,\,\, \phi(x,0)=f(x) +
\displaystyle \sum_{i=1}^N\rho_i^2\left[\xi_i(x,T)\right]~~~~\mbox{in}~~~~
\Omega.
\end{array}\right.
\end{equation}
Then, the following result holds:
\begin{theorem} \label{teo3.1} For each $1,2, \ldots, N,$ assume that  $\rho_i\in L^{\infty}(\Omega)$, $f\in H$, and that  $\sum_{i=1}^N\alpha_i$ is sufficiently small. Let $g:[0,\infty)\rightarrow[0,\infty)$ be given by $g(t)=\gamma e^{-\delta t},$ where $\gamma=\delta(\nu-\mu),$  $\mu=k\lambda^{-1},$ and
	$\delta=\lambda^{-1}$ are positive constants. Then, system $(\ref{eqA.1})$ admits a unique solution $$\{\phi,\xi_i,\widehat{p},p\} \in C^0([0,T];H)\times
	C^0([0,T];V)\times L^2(0,T;L^2(\Omega))\times  L^2(0,T;H^1(\Omega))$$ such that
\begin{equation*}
	\begin{array}{c}
\displaystyle \displaystyle \frac{\partial \phi}{\partial
		t} - \mu\Delta \phi - \displaystyle \int_t^Tg(\eta-t)\Delta \phi(\eta)d\eta + \nabla \widehat{p}=0\,\,\,
	\mbox{in}\,\,\,\,\, L^2(0,T;\textbf{L}^2(\Omega)),\vspace{0.3cm}\\
	\displaystyle \frac{\partial \xi_i}{\partial
	t}-\mu\Delta\xi_i - \int_0^tg(t-\sigma)\Delta
\xi_i(\sigma)d\sigma +  \nabla {p}=-\alpha_i\phi\chi_{\mathcal{O}_i} \,\,\mbox{in}\,\,L^2(0,T;\textbf{L}^2(\Omega)),\vspace{0.3cm}\\
	\xi_i(0)=0\,\,\,\,\,\mbox{and}\,\,\,\,\,\phi(0)=f+
	\displaystyle \sum_{i=1}^N\rho_i^2[\xi_i(T)]\,\,\,\,\mbox{in}\,\,\,\,\Omega.
	\end{array}
\end{equation*}
\end{theorem}

\proof  We proceed by using the Faedo-Galerkin method for a special basis
	$\{w_j\}_{j\in \mathbb{N}}$ of $V$. We refer to the book of Temam (1979) for a complete description of this method. Let $V_m=[w_1,....,w_m]$  be the subspace generated by the first $m$ vectors $\{w_j\}_j.$ The approximate problem consists in finding functions $\phi_m(x,t)=\sum_{j=1}^mq_{jm}(t)w_j(x),$  and $\xi_{im}(x,t)=\sum_{j=1}^mh^i_{j m}(t)w_j(x),$ where  $q_{jm}(t)$ and $h^i_{j m}(t)$ are real functions defined in $[0,T],$ such that:
\begin{equation}\label{eqA.2}
	\left|\begin{array}{l}
\left  ( \displaystyle \frac{\partial \phi_m(t)}{\partial
	t},w_j\right) + \mu\left((\phi_m(t),w_j)\right)\;   +
	\displaystyle\int_t^Tg(\eta-t)\left((\phi_m(\eta),w_j)\right)d\eta =0,\vspace{0.3cm}\\
	\left  ( \displaystyle \frac{\partial \xi_{im}(t)}{\partial
	t},w_j\right) + \mu\left((\xi_{im}(t),w_j)\right)+
	\displaystyle\int_0^tg(t-\sigma)\left((\xi_{im}(\sigma),w_j)\right)d\sigma \\
	=-\alpha_i\left(\phi_m(t)\chi_{\mathcal{O}_i},w_j\right),\\
	\displaystyle\xi_{im}(0)=0,\,\,\phi_{m}(0) -
	\sum_{i=1}^N\rho_i^2\xi_{im}(T)=\sum_{\lambda=1}^m\left(f,w_j\right)w_j,\,\mbox{for
		all}\,w_j\in V_m.
	\end{array}\right.
	\end{equation}
	Multiplying $(\ref{eqA.2})_1$ by $q_{jm}(t)$ and summing for $j=1,.....m$,
	we find that:
	\begin{equation}\label{eqA.6}
	\begin{array}{c}
	\begin{array}{c}
	\displaystyle\frac{d}{dt}|\phi_m(t)|^2 + \mu||\phi_m(t)||^2 +
	\int_0^tg(t - \sigma)((\phi_m(\sigma),\phi_m(t)))d\sigma\leq\vspace{0.3cm}
	\displaystyle|\phi_m(t)|^2.
	\end{array}
	\end{array}
	\end{equation}
	Integrating $(\ref{eqA.6})$ from $0$ to $t,$ and applying Fubini's formula and
	Gronwall's inequality, we get that:
	\begin{equation}\label{eqA.7}
	\begin{array}{l}
	|\phi_m(t)|^2\leq C(T)|\phi_m(0)|^2,\,\,\mbox{for all}\,\,m,\,t\geq
	0.
	\end{array}
	\end{equation}
	Analogously, multiplying $(\ref{eqA.2})_2$ by $\lambda_jh^i_{jm}(t)$ and summing from $j=1,......,m$, we obtain:
	\begin{equation}\label{eqA.8}
	\begin{array}{l}
	\displaystyle \frac{d}{dt}||\xi_{im}(t)||^2 + \frac{\mu}{2}|\Delta
	\xi_{im}(t)|^2  +
	\int_0^tg(t-\sigma)(\Delta\xi_m(\sigma),\Delta\xi_m(t))d\sigma\leq
	\displaystyle |\phi_m(t)|^2.
	\end{array}
	\end{equation}
	Thus, combining $(\ref{eqA.7}),$ and $(\ref{eqA.8})$ it follows that:
	\begin{equation}\label{eqA.9}
	\begin{array}{l}
	\displaystyle \frac{d}{dt}||\xi_{im}(t)||^2 + \frac{\mu}{2}|\Delta
	\xi_{im}(t)|^2 +
	\int_0^tg(t-\sigma)(\Delta\xi_m(\sigma),\Delta\xi_m(t))d\sigma\leq
	C(T)|\phi_m(0)|^2.
	\end{array}
	\end{equation}
	Integrating $(\ref{eqA.9})$ from $0$ to $t$, and combining again Fubini's formula and
	Gronwall's inequality, we have that:
	\begin{equation}\label{eqA.10}
	\begin{array}{l}
	\displaystyle\|\xi_{im}(t)\|^2\leq
	\tilde{C}(T)|\phi_m(0)|^2\,\,\mbox{for all}\,\,m,\,t\in[0,T].
	\end{array}
	\end{equation}
	In particular,
	\begin{equation}\label{eqA.11}
	\begin{array}{l}
	\displaystyle\|\xi_{im}(T)\|^2\leq
	\tilde{C}(T)|\phi_m(0)|^2\,\,\mbox{for all}\,\,m.
	\end{array}
	\end{equation}

	From the definition of $\phi_m(0),$ (see  $(\ref{eqA.2})$), we see that:
	\begin{equation}\label{eqA.12}
	\begin{array}{lll}
	\displaystyle|\phi_m(0)|^2 &\leq&  c_0C_1(N)\displaystyle\left(\sum_{i=1}^N|\rho_i|^4_{L^{\infty}(\Omega)}\right)\left(\sum_{i=1}^N\|\xi_{im}(T)\|^2\right)+ |f|^2,
	\end{array}
	\end{equation}
	where $c_0$ is the constant of immersion of the  $V$ into  $H.$

	From $(\ref{eqA.11})$ and $(\ref{eqA.12}),$ we conclude that:
	\begin{equation}\label{eqA.13}
	\begin{array}{l}
	\displaystyle
	\left[ 1-c_0C_1(N)\displaystyle \left(\sum_{i=1}^N|\rho_i|^4_{L^{\infty}(\Omega)}\right)
	\sum_{i=1}^N\frac{\alpha_i^2}{\mu} \right]\sum_{i=1}^N\|\xi_{im}(T)\|^2\leq
	\sum_{i=1}^N\frac{\alpha_i^2} {\mu}|f|^2.
	\end{array}
	\end{equation}
	
	Notice that we can assume $\sum_{i=1}^N\alpha_i^2$ sufficiently small so that:
	$$\beta=1-c_0C_1(N)
	\displaystyle\left(\sum_{i=1}^N|\rho_i|^4_{L^{\infty}(\Omega)}\right)\sum_{i=1}^N\frac{\alpha_i^2}
	{\mu}>0.$$

Hence, from (\ref{eqA.13}), we obtain that:
	\begin{equation}\label{eqA.14}
	\begin{array}{l}
	\displaystyle \sum_{i=1}^N\|\xi_{im}(T)\|^2\leq
	\frac{1}{\beta\mu}\sum_{i=1}^N\alpha_i^2|f|^2.
	\end{array}
	\end{equation}

	Now, combining  $(\ref{eqA.12})$ and $(\ref{eqA.14}),$  we deduce that:
	\begin{equation}\label{eqA.15}
	\begin{array}{c}
	\displaystyle|\phi_m(0)|^2\leq \displaystyle
	\left[c_0C_1(N)\displaystyle\left(\sum_{i=1}^N|\rho_i|^4_{L^{\infty}(\Omega)}\right)
	\frac{1}{\beta\mu}\sum_{i=1}^N\alpha_i^2 +
	1\right]|f|^2,\,\,\mbox{for all}\,\,m.
	\end{array}
	\end{equation}

    From  $(\ref{eqA.6})-(\ref{eqA.10}),$ it follows that:
	\begin{equation}\label{eqA.16}
	\begin{array}{l}
	\left(\phi_m\right)_{m\in \mathbb{N}}\,\,\,\mbox{is bounded in}\,\,\,L^2(0,T;V),\\
	\left(\phi_m\right)_{m\in \mathbb{N}}\,\,\,\mbox{is bounded in}\,\,\,L^{\infty}(0,T;H),\\
	\left(\xi_{im}\right)_{m\in \mathbb{N}}\,\,\,\mbox{is bounded in}\,\,\,L^{\infty}(0,T;V),\\
	\left(\Delta\xi_{im}\right)_{m\in \mathbb{N}}\,\,\,\mbox{is bounded
		in}\,\,\,L^2(0,T;H),
	\end{array}
	\end{equation}
	for $i=1,2,3,.....,N.$ In particular, from $(\ref{eqA.16})_4,$ we conclude that:
	\begin{equation}\label{eqA.17}
	\begin{array}{l}
	\left(\xi_{im}\right)_{m\in \mathbb{N}}\,\,\,\mbox{is bounded
		in}\,\,\,L^2(0,T;\textbf{H}^2(\Omega)).
	\end{array}
	\end{equation}
	
	Moreover, from $(\ref{eqA.16})_1-(\ref{eqA.16})_3,$ and $(\ref{eqA.17}),$  we can extract
	subsequences of $(\phi_m)_{m\in \mathbb{N}}$ and $(\xi_{im})_{m\in \mathbb{N}}$ such that:
	\begin{equation}\label{eqA.18}
	\begin{array}{l}
	\phi_m\rightharpoonup \phi\,\,\,\mbox{in}\,\,\,L^{2}(0,T;V),\\
	\phi_m\stackrel{*}{\rightharpoonup} \phi\,\,\,\mbox{in}\,\,\,L^{\infty}(0,T;H),\\
	\xi_{im}\stackrel{*}{\rightharpoonup}
	\xi_i\,\,\,\mbox{in}\,\,\,L^{\infty}(0,T;V),\\
	\xi_{im}\rightharpoonup
	\xi_i\,\,\,\mbox{in}\,\,\,L^{2}(0,T;\textbf{H}^2(\Omega)),\\
	\end{array}
	\end{equation}
	for $i=1,2,....,N.$
	
	From $(\ref{eqA.18})_1$ and $(\ref{eqA.18})_2$, we find that:
	\begin{equation*}
	\begin{array}{c}
	\displaystyle\frac{d}{dt}(\phi(t),v)\;  + \mu(( \phi(t),v)) +
	\displaystyle\int_0^tg(t - \sigma)\left((\phi(\sigma),v)\right)d\sigma =
	0\,\,\,\,\, \mbox{for all}\,\,\,\ \,v\,\,\in V\,\,\,\mbox{in}\,\,\,\mathcal{{D}}'(0,T).
	\end{array}
	\end{equation*}

This means that:
	\begin{equation}\label{eqA.19}
	\begin{array}{c}
	\displaystyle \frac{\partial \phi}{\partial
		t}\in L^2(0,T;V'),\vspace{0.3cm}\\
	\displaystyle \frac{\partial \phi}{\partial
		t}\ -\mu \Delta \phi  -\displaystyle \int_t^Tg(\eta-t)\Delta \phi(\eta)d\eta  =0\vspace{0.3cm}\,\,\, \mbox{in the
		sense of }\,\,\,\,\,\,L^2(0,T;V').
	\end{array}
	\end{equation}

	Moreover,
	\begin{equation}\label{eqA.20}
	\begin{array}{l}
	 \displaystyle \frac{\partial \phi}{\partial
		t} - \mu \Delta \phi(t) -\int_t^Tg(\eta-t)\Delta \phi(\eta)d\eta  =0,
	\end{array}
	\end{equation}
	in the sense of $V'$ in $[0,T]$. In particular, we have that:
	\begin{equation}\label{eqA.21}
	\begin{array}{l}
	\displaystyle \left< \displaystyle \frac{\partial \phi}{\partial
		t} - \mu \Delta \phi(t) -\int_t^Tg(\eta-t)\Delta
	\phi(\eta)\,d\eta, v \right >  =0,
	\end{array}
	\end{equation}
	for all $v\in \mathcal{V}$ in $[0,T]$. Hence, as in Temam (1979), there
	exists $\widehat{p}(t)\in L^2(\Omega)$ such that:
	\begin{equation}\label{eqA.22}
	|\widehat{p}(t)|_{L^2(\Omega)}\leq C ||\nabla
	\widehat{p}(t)||_{\textbf{H}^{-1}(\Omega)},
	\end{equation}
	and
	\begin{equation}\label{eqA.23}
	-\nabla \widehat{p}(t)= \displaystyle \frac{\partial \phi}{\partial
		t} - \mu \Delta \phi(t) -  \int_t^Tg(\eta-t)\Delta \phi(\eta)\,d\eta.
	\end{equation}
	Therefore, from  $(\ref{eqA.22})$ and  $(\ref{eqA.23}),$ we conclude that:
	\begin{equation}\label{eqA.24}
	\nabla \widehat{p} \in L^2(0,T;\textbf{H}^{-1}(\Omega)),
	\end{equation}
	and
	\begin{equation}\label{eqA.25}
	\widehat{p}\in L^2(0,T;L^2(\Omega)).
	\end{equation}
	Since $ \displaystyle \frac{\partial \phi}{\partial
		t}\in L^2(0,T;V')$, $\phi\in L^2(0,T;V),$ and $V\subset H\subset V'$, then  as in Temam (1979),  we see that:
	\begin{equation}\label{eqA.26}
	\phi\in C^0([0,T]; H).
	\end{equation}
	
Finally, we deduce  $(\ref{eqA.18})$ that:
	\begin{equation}\label{eqA.27}
	\xi_i\in L^2(0,T;\textbf{L}^2(\Omega)),
	\end{equation}
	\begin{equation}\label{eqA.28}
	p\in L^2(0,T;H^1(\Omega)),
	\end{equation}
	\begin{equation}\label{eqA.29}
	\xi_i\in C^0([0,T];V),
	\end{equation}
	\begin{equation}\label{eqA.30}
	\begin{array}{c}
	\displaystyle \displaystyle \frac{\partial \phi}{\partial
		t} - \mu\Delta \phi - \displaystyle \int_t^Tg(\eta-t)\Delta \phi(\eta)d\eta + \nabla \widehat{p}=0\,\,\,
	\mbox{in}\,\,\,\,\, L^2(0,T;\textbf{L}^2(\Omega)),
	\end{array}
	\end{equation}
	\begin{equation}\label{eqA.31}
	\begin{array}{c}
	\displaystyle \displaystyle \frac{\partial \xi_i}{\partial
		t} - \mu\Delta \xi_i -\int_0^tg(t-\sigma)\Delta \xi_i(\sigma)d\sigma + \nabla
	p= - \alpha_i\chi_{\mathcal{O}_i}\,\,\, \mbox{in}\,\,\,\,\,
	L^2(0,T;\textbf{L}^2(\Omega)),
	\end{array}
	\end{equation}
	for all $i=1,2,....,N.$ The uniqueness is obtained of standard form. \cqdf

\begin{remark}
The  proof of the existence, uniqueness, and regularity of solutions for the system $(\ref{eqI.9})$ follows similarly to the demonstration of the previous theorem, with some suitable adaptations.
Due to that, we will omit this proof.
\end{remark}

\paragraph{\bf Acknowledgement} The authors want to express their gratitude to the anonymous reviewers for their questions and commentaries; they were
very helpful in improving this article. Isa\'{\i}as Pereira  de Jesus was supported by CNPq/Brazil grants 307488/2019-5, 305394/2022-3 and PRPG (UFPI).

\paragraph

\end{document}